\documentclass[11pt]{amsart}
\usepackage{amsmath,amsthm}
\newtheorem{theorem}{Theorem}[section]
\newtheorem{cor}[theorem]{Corollary}
\newtheorem{lemma}[theorem]{Lemma}
\newtheorem{prop}[theorem]{Proposition}

\theoremstyle{remark}
\newtheorem{remark}[theorem]{Remark}

\theoremstyle{definition}
\newtheorem{definition}[theorem]{Definition}
\newtheorem{notation}[theorem]{Notation}
\numberwithin{equation}{section}
\DeclareMathOperator{\Ind}{Ind}
\DeclareMathOperator{\dashind}{-Ind}
\DeclareMathOperator{\End}{End}
\DeclareMathOperator{\cv}{cov}
\DeclareMathOperator{\ran}{ran}
\DeclareMathOperator{\id}{id}
\DeclareMathOperator{\Span}{span}
\DeclareMathOperator{\clsp}{\overline{span}}
\newsymbol\rtimes 216F
\newcommand{\norm}[1]{\lVert#1\rVert}
\newcommand{\cp}{C\textstyle{\cross_{\beta,X}} P}
\newcommand{\bpp}{B_P \textstyle{\cross_{\tau, X}} P}
\newcommand{\Star}{${}^*$\ndash}
\newcommand{\cross}{\rtimes}
\newcommand{\ndash}{\nobreakdash-}
\newcommand{\field}[1]{\mathbb{#1}}
\newcommand{\CC}{\field{C}}
\newcommand{\FF}{\field{F}}
\newcommand{\NN}{\field{N}}

\newcommand{\TT}{\field{T}}
\newcommand{\ZZ}{\field{Z}}
\newcommand{\Aa}{{\mathcal A}}
\newcommand{\Bb}{{\mathcal B}}

\newcommand{\Ff}{{\mathcal F}}
\newcommand{\Gg}{{\mathcal G}}
\newcommand{\Hh}{{\mathcal H}}
\newcommand{\Ii}{{\mathcal I}}
\newcommand{\Kk}{{\mathcal K}}
\newcommand{\Ll}{{\mathcal L}}

\newcommand{\Oo}{{\mathcal O}}
\newcommand{\Pp}{{\mathcal P}}
\newcommand{\Tt}{{\mathcal T}}
\newcommand{\Uu}{{\mathcal U}}

\begin{document}
%
%
\title[Product systems of Hilbert bimodules]
{Discrete product systems of Hilbert bimodules}
\author[Neal J. Fowler]{Neal J. Fowler}
\address{Department of Mathematics  \\
      University of Newcastle\\  NSW  2308\\ AUSTRALIA}
\email{neal@math.newcastle.edu.au}
\date{April 19, 1999}
\thanks{The author thanks the University of Victoria, Canada, for its
hospitality while this research was being completed.}

\begin{abstract}
A Hilbert bimodule is a right Hilbert module $X$ over a $C^*$-algebra $A$
together with a left action of $A$ as adjointable operators on $X$.
We consider families $X = \{X_s :s\in P\}$ of Hilbert bimodules,
indexed by a semigroup $P$, which are endowed with a multiplication
which implements isomorphisms $X_s\otimes_A X_t \to X_{st}$;
such a family is a called a product system.
We define a generalized Cuntz-Pimsner algebra $\Oo_X$,
and we show that every twisted crossed product
of $A$ by $P$ can be realized as $\Oo_X$
for a suitable product system $X$.
Assuming $P$ is quasi-lattice ordered in the sense of Nica,
we analyze a certain Toeplitz extension $\Tt_{\cv}(X)$ of $\Oo_X$
by embedding it in a crossed product $\bpp$
which has been ``twisted'' by $X$; our main Theorem is a characterization
of the faithful representations of $\bpp$.
\end{abstract}

\maketitle
%
%

\section*{Introduction}

Suppose $X$ is a right Hilbert module over a $C^*$\ndash algebra $A$.
If $X$ also carries a left action of $A$ as adjointable operators on $X_A$,
we  call $X$ a {\em Hilbert bimodule\/} over $A$.
In \cite{pimsner}, Pimsner associated with every such bimodule $X$ a $C^*$\ndash
algebra $\Oo_X$, which we shall call the {\em Cuntz-Pimsner algebra\/} of $X$,
and showed that every crossed product by $\ZZ$
and every Cuntz-Krieger algebra can be realized as $\Oo_X$ for suitable $X$.
He also commented that the algebras $\Oo_X$ include the crossed products by $\NN$;
that is, for each endomorphism $\alpha$ of a $C^*$\ndash algebra
$A$ there is a bimodule $X = X(\alpha)$ such that $\Oo_X$ is
canonically isomorphic to the semigroup crossed product
$A\cross_\alpha\NN$ of \cite{cuntz82,stacey}.

The work in this paper is motivated by the following observation,
which also serves as our primary example.
Suppose $\beta$ is an action of a discrete semigroup $P$
as endomorphisms of a $C^*$\ndash algebra $A$.
For each $s\in P$ let $X_s := X(\beta_s)$ be the bimodule canonically associated
with the endomorphism $\beta_s$.
Then the family $X = \{X_s: s\in P\}$ admits an associative  multiplication
$(x,y)\in X_s\times X_t \mapsto xy\in X_{ts}$ which implements
isomorphisms $X_s\otimes_A X_t \to X_{ts}$;
we call a family with this structure a {\em product system of Hilbert bimodules\/}.
(In this example $X$ is a product system over the opposite semigroup $P^o$.)
Such families generalize the product systems of
\cite{dinhjfa,dinhjot,fowrae,fowler},
where the fibers $X_s$ are complex Hilbert spaces (bimodules over $\CC$).

To each product system $X$ we associate a generalized Cuntz-Pimsner
algebra $\Oo_X$.  When $X$ is the product system associated with
the semigroup dynamical system $(A,P,\beta)$, $\Oo_X$ is canonically
isomorphic to the semigroup crossed product $A\cross_\beta P$.
Moreover, if we ``twist'' $X$ by a multiplier $\omega:P\times P\to\TT$,
then the corresponding Cuntz-Pimsner algebra is isomorphic to
the twisted semigroup crossed product $A\cross_{\beta,\omega} P$.
Our construction applies even when $A$ is nonunital provided each
endomorphism $\beta_s$ extends to the multiplier algebra $M(A)$.

The aim of this paper is to take a first step towards analyzing
the Cuntz-Pimsner algebra of a product system $X$.
Following Pimsner \cite{pimsner},
we begin by studying the structure of its Toeplitz extension $\Tt_X$.
This algebra is universal for {\em Toeplitz  representations\/} of $X$;
these are multiplicative maps whose restriction to each fiber $X_s$
is a Toeplitz representation in the sense of \cite{fowrae2}.
Our results generalize those of \cite{fowrae} for product systems
of Hilbert spaces; indeed, much of the paper is devoted to
adapting the methods of \cite{fowrae} to the bimodule setting.
Thus our basic assumptions about the underlying semigroup $P$ are
as in \cite{fowrae}:
to allow our analysis to extend beyond the totally-ordered case,
we assume that $P$ is the positive cone of a group $G$
such that $(G,P)$ is quasi-lattice ordered in the sense of Nica \cite{nica}.
The class of such $(G,P)$ includes all direct sums and free products
of totally ordered groups.
We also impose a covariance condition,
called {\em Nica covariance\/}, on Toeplitz representations of $X$.
This means that the universal $C^*$\ndash algebra $\Tt_{\cv}(X)$
which we analyze is in general a quotient of $\Tt_X$.
However, if $(G,P)$ is totally-ordered,
then Nica-covariance is automatic,
and hence $\Tt_{\cv}(X)$ is the same as $\Tt_X$.

Our main goal is to characterize the faithful representations
of $\Tt_{\cv}(X)$.
We accomplish this by embedding $\Tt_{\cv}(X)$
in a certain twisted semigroup crossed product $\bpp$
(Theorem~\ref{theorem:subalgebra}), and then characterizing
its faithful representations
(Theorem~\ref{theorem:faithfulness of representations}).
When $P = \NN$, $\Tt_{\cv}(X)$ is precisely the Toeplitz algebra of
the Hilbert bimodule $X_1$ (the fiber over $1\in\NN$),
and our Theorem~\ref{theorem:faithfulness of representations}
reduces to \cite[Theorem~2.1]{fowrae2}.
In fact, the analysis in \cite{fowrae2} was motivated
by our preliminary work on this paper.
We would like to point out in particular how
the stronger result \cite[Theorem~3.1]{fowrae2}
arose from our investigations into product systems,
for it serves as a good illustration of the usefulness
of Nica covariance.
Suppose $Z$ is an orthogonal direct sum
$\bigoplus_{\lambda\in\Lambda} Z^\lambda$ of Hilbert bimodules.
Let $G$ be the free group on $\Lambda$,
let $P$ be the subsemigroup of $G$ generated by $\Lambda$,
and let $X$ be the unique product system over $P$
whose fiber over $\lambda$ is $Z^\lambda$.
Then $\Tt_{\cv}(X)$ is canonically isomorphic
to the Toeplitz algebra of the bimodule $Z$,
and \cite[Theorem~3.1]{fowrae2} follows from our
Theorem~\ref{theorem:faithfulness of representations}.

The main application of \cite[Theorem~3.1]{fowrae2}
was to establish the simplicity of the graph algebras
associated with certain infinite directed graphs
\cite[Corollary~4.3]{fowrae2}.
Although here we confine our applications
to twisted semigroup crossed products,
we anticipate that our results will also give interesting
information about $\Oo_X$ when each of the fibers of $X$
arise from infinite directed graphs. 

We begin in Section~\ref{section:prelims} by giving a
brief review of Hilbert bimodules,
their representations, and their $C^*$\ndash algebras.
In Section~\ref{section:product systems} we introduce
product systems of Hilbert bimodules,
discuss their representations, and define
the algebras $\Tt_X$ and $\Oo_X$.
In Section~\ref{section:crossed products} we associate with each twisted
semigroup dynamical system $(A,P,\beta,\omega)$ a product system
$X = X(A,P,\beta,\omega)$ whose Cuntz-Pimsner algebra $\Oo_X$
is the twisted semigroup crossed product $A\cross_{\beta,\omega} P$.
We show that the Toeplitz algebra of $X(A,P,\beta,\omega)$
also has a crossed product structure, and this motivates the definition
of a ``Toeplitz'' crossed product $\Tt(A\cross_{\beta,\omega} P)$
in which the endomorphisms are implemented not by isometries,
but rather by partial isometries.

In Section~\ref{section:crossed products2} we generalize the notion
of twisted crossed product by replacing the multiplier
$\omega$ by a product system $X$ of Hilbert bimodules.
This extends the philosophy developed in \cite{fowrae}
that one should regard product systems as noncommutative cocycles.
Hence given an action $\beta$ of $P$ as endomorphisms of a $C^*$\ndash
algebra $C$, we consider $(C,P,\beta,X)$ as a twisted
semigroup dynamical system,
and we define a twisted crossed product $C\cross_{\beta,X} P$.

In Section~\ref{section:covariance} we assume that $(G,P)$
is quasi-lattice ordered, and we discuss the notion
of Nica covariance for a Toeplitz representation.
As illustrated in \cite[Example~1.3]{fowler}
using product systems of Hilbert spaces,
when $(G,P)$ is not a total order
it is possible that the $C^*$\ndash algebra $\Tt_{\cv}(X)$
which is ``universal'' for such representations
may admit representations which are not the integrated form
of a Nica-covariant Toeplitz representation.
To avoid this pathology
we adapt the methods of \cite{fowler} to our setting:
we define the notion of a product system being
{\em compactly aligned\/}, and show that $\Tt_{\cv}(X)$
is truly universal when $X$ is compactly aligned
(Proposition~\ref{prop:preserve covariance}).
We show that $X$ is compactly aligned
if the left action of $A$ on each fiber $X_s$ is by compact
operators (Proposition~\ref{prop:CA conditions});
it follows that the product systems $X(A,P,\beta,\omega)$
associated with twisted semigroup dynamical systems
are compactly aligned.

In Section~\ref{section:system} we consider a certain
$C^*$\ndash subalgebra $B_P$ of $\ell^\infty(P)$
which is invariant under left translation $\tau:P\to\End(\ell^\infty(P))$.
As in \cite{lacarae,fowrae}, covariant representations
of the twisted system $(B_P,P,\tau,X)$
are in one-one correspondence with Toeplitz representations of $X$
which are Nica-covariant (Proposition~\ref{prop:Lpsi}),
and hence $\Tt_{\cv}(X)$ embeds naturally as a subalgebra
of $\bpp$ (Theorem~\ref{theorem:subalgebra}).
When the left action of $A$ on each fiber $X_s$ is by compact operators,
$\Tt_{\cv}(X)$ is all of $\bpp$.

In Section~\ref{section:faithful} we prove our main result,
Theorem~\ref{theorem:faithfulness of representations},
which characterizes the faithful
representations of $B_P\cross_{\tau,X} P$
under the assumption that $X$ is compactly aligned
and $(B_P,P,\tau,X)$ satisfies a certain amenability hypothesis.
In Section~\ref{section:amenability} we give conditions on $(G,P)$
which ensure that $(B_P,P,\tau,X)$ is amenable.
In particular, $(B_P,P,\tau,X)$ is amenable
if $X$ is compactly aligned and
$(G,P)$ is a free product $*(G^\lambda,P^\lambda)$
with each $G^\lambda$ amenable (Corollary~\ref{cor:amenability}).

Finally, in Section~\ref{section:applications} we apply our
Theorem~\ref{theorem:faithfulness of representations}
to the product system
$X(A,P,\beta,\omega)$ of Section~\ref{section:crossed products}.
When $(G,P)$ is a total order,
$\bpp$ is isomorphic to the Toeplitz crossed product
$\Tt(A\cross_{\beta,\omega} P)$;
in general $\bpp$ is a certain quotient $\Tt_{\cv}(A\cross_{\beta,\omega} P)$
which also has a crossed product structure,
and Theorem~\ref{theorem:crossed products}
characterizes its faithful representations.
Applying this to the twisted system $(B_P,P,\tau,\omega)$,
we show that $\Tt_{\cv}(B_P\cross_{\tau,\omega} P)$ is universal
for partial isometric representations of $P$ which are {\em bicovariant\/}
(Proposition~\ref{prop:universal bicovariant}),
and we obtain a characterization of its faithful representations
(Theorem~\ref{theorem:bicovariant})
which is particularly nice
when $P$ is the free semigroup on infinitely many generators
(Theorem~\ref{theorem:free}).

The author thanks Iain Raeburn for the many helpful discussions
while this research was being conducted.

\section{Preliminaries}\label{section:prelims}

Let $A$ be a separable $C^*$\ndash al\-ge\-bra.
A {\em Hilbert bimodule over $A$\/} is a right Hilbert $A$-module $X$
together with a \Star homomorphism $\phi:A\to\Ll(X)$
which is used to define a left action of $A$ on $X$ via $a\cdot x := \phi(a)x$
for $a\in A$ and $x\in X$.
A {\em Toeplitz representation} of $X$
in a $C^*$\ndash al\-ge\-bra $B$ is a pair $(\psi,\pi)$
consisting of a linear map $\psi:X\to B$
and a homomorphism $\pi:A\to B$ such that
\begin{align*}
\psi(x\cdot a) &= \psi(x)\pi(a), \\
\psi(x)^*\psi(y)&= \pi(\langle x,y \rangle_A),\ \text{ and} \\
\psi(a\cdot x) &= \pi(a)\psi(x)
\end{align*}
for $x,y\in X$ and $a\in A$.
Given such a representation, there is
homomorphism $\pi^{(1)}:\Kk(X)\to B$ which satisfies
\begin{equation}\label{eq:pione}
\pi^{(1)}(\Theta_{x,y}) = \psi(x)\psi(y)^*
\qquad\text{for all $x,y\in X$,}
\end{equation}
where $\Theta_{x,y}(z) := x\cdot\langle y,z \rangle_A$ for $z\in X$;
see \cite[p.~202]{pimsner},\cite[Lemma~2.2]{kpw}, and \cite[Remark~1.7]{fowrae2}
for details.
We say that $(\psi,\pi)$ is {\em Cuntz-Pimsner covariant\/}
if
\[
\pi^{(1)}(\phi(a)) = \pi(a)
\qquad\text{for all $a \in \phi^{-1}(\Kk(X))$.}
\]
The {\em Toeplitz algebra\/} of $X$ is the $C^*$\ndash al\-ge\-bra $\Tt_X$
which is universal for Toeplitz representations of $X$
\cite{pimsner,fowrae2}, and the {\em Cuntz-Pimsner algebra\/}
of $X$ is the $C^*$\ndash al\-ge\-bra $\Oo_X$ which is universal for
Toeplitz representations which are Cuntz-Pimsner covariant
\cite{pimsner,dpz,kpw,ms, ms2,fmr}.

Every right Hilbert $A$-module $X$ is essential, in the sense that
$X$ is the closed linear span of elements $x\cdot a$.
We say that a Hilbert bimodule $X$ is {\em essential\/}
if it is also essential as a left $A$-module; that is, if
\[
X = \clsp\{\phi(a)x: a\in A, x\in X\}.
\]
When $X$ is essential,
two applications  of the Hewitt-Cohen Factorization Theorem allow us to write
any $x\in X$ as $\phi(a)y\cdot b$  for some $y\in X$ and $a,b\in A$.
Hence if $(a_i)$ is an approximate identity for  $A$,
then
\begin{equation}\label{eq:approximate identity}
\lVert x - x\cdot a_i\rVert \to 0
\quad\text{and}\quad
\lVert x - \phi(a_i)x\rVert \to 0
\qquad\text{for all $x\in X$.}
\end{equation}

\section{Product systems of Hilbert bimodules}
\label{section:product systems}

For each $n\ge 1$ the $n$-fold internal tensor product
$X^{\otimes n} := X\otimes_A \dotsb \otimes_A X$ has a natural structure
as a Hilbert bimodule over $A$; see \cite[ Section~2.2]{ms2} for details.
The following definition,
based on Arveson's continuous tensor product systems over $(0,\infty)$ \cite{arv},
generalizes the collection $\{X^{\otimes n}: n\in\NN\}$
to semigroups other than $\NN$.

\begin{definition}\label{defn:ps}
Suppose $P$ is a countable semigroup with identity $e$
and $p:X\to P$ is a family of Hilbert bimodules over $A$.
Write $X_s$ for the fibre $p^{-1}(s)$ over $s\in P$,
and write $\phi_s:A\to\Ll(X_s)$ for the homomorphism which
defines the left action of $A$ on $X_s$.
We say that $X$ is a {\em (discrete) product system over $P$\/} if $X$ is a semigroup,
$p$ is a semigroup homomorphism,
and for each $s,t\in P\setminus\{e\}$ the map
$(x,y) \in X_s \times X_t \mapsto xy \in X_{st}$
extends to an isomorphism of the Hilbert bimodules
$X_s\otimes_A X_t$ and $X_{st}$.
We also require that $X_e = A$
(with its usual right Hilbert module structure and $\phi_e(a)b = ab$ for $a,b\in A$),
and that the multiplications $X_e\times X_s\to X_s$
and $X_s\times X_e \to X_s$ satisfy
\begin{equation}\label{eq:multiply}
ax = \phi_s(a)x,\qquad xa = x\cdot a
\qquad \text{for $a\in X_e$ and $x\in X_s$.}
\end{equation}
\end{definition}

\begin{remark}\label{remark:essential}
Multiplication $X_e \times X_s \to X_s$ will not induce an isomorphism
$X_e\otimes_A X_s \to X_s$ unless $X_s$ is essential as a left $A$-module.
\end{remark}

\begin{remark} 
The associativity of multiplication in $X$ implies that
$\phi_{st}(a) = \phi_s(a)\otimes_A 1^t$  for all $a\in A$;
that is, $\phi_{st}(a)(xy) = (\phi_s(a)x)y$
for all $x\in X_s$ and $y\in X_t$.
\end{remark}

\begin{remark} It is possible that some of the $X_s$ may be zero.
\end{remark}

\begin{definition}\label{defn:psrep}
Suppose $B$ is a $C^*$\ndash al\-ge\-bra and $\psi:X\to B$;
write $\psi_s$ for the restriction of $\psi$ to $X_s$.
We call $\psi$ a {\em Toeplitz representation\/} of $X$ if

\textup{(1)} For  each $s\in P$, $(\psi_s, \psi_e)$
is a Toeplitz representation of $X_s$, and

\textup{(2)} $\psi(xy) = \psi(x)\psi(y)$ for $x,y\in X$.

\noindent If in addition each $(\psi_s,\psi_e)$ is Cuntz-Pimsner covariant,
we say that $\psi$ is {\em Cuntz-Pimsner covariant\/}.
\end{definition}

\begin{remark}\label{remark:toeplitz rep}
By \cite[Remark~1.1]{fowrae2}, every Toeplitz representation $\psi$ is contractive;
moreover, if the homomorphism $\psi_e:A\to B$ is isometric, then so is $\psi$.
Also, since we are assuming \eqref{eq:multiply}, a map $\psi:X\to B$
is a Toeplitz representation if it satisfies both (2) and

\textup{(1')} $\psi_s(x)^*\psi_s(y) = \psi_e(\langle x,y \rangle_A)$
whenever $s\in P$ and $x,y\in X_s$.
\end{remark}

\begin{notation}
We write $\psi^{(s)}$ for the homomorphism of $\Kk(X_s)$ into $B$
which corresponds to the pair $(\psi_s,\psi_e)$,
as in \eqref{eq:pione}; that is,
\[
\psi^{(s)}(\Theta_{x,y}) = \psi_s(x)\psi_s(y)^*
\qquad\text{for all $x,y\in X_s$.}
\]
\end{notation}

\subsection*{The Fock representation}
Let $F(X)$ be the right Hilbert $A$-module
\[
F(X) := \bigoplus_{s\in P} X_s.
\]
By this we mean the following:  as a set,
$F(X)$ is the subset of $\prod_{s\in P} X_s$
consisting of all elements $(x_s)$ for which
$\sum_{s\in P} \langle x_s, x_s \rangle_A$ is summable in $A$;
that is, for which $\sum_{s\in F} \langle x_s, x_s \rangle_A$
converges in norm as $F$ increases over the finite subsets of $P$.
We write $\oplus x_s$ for $(x_s)$ to indicate that the above series is summable.
The right action of $A$ is given by $(\oplus x_s)\cdot a := \oplus (x_s \cdot a)$,
and the inner product by
$\langle \oplus x_s, \oplus y_s \rangle_A
 := \sum_{s\in P} \langle x_s, y_s \rangle_A$.
The algebraic direct sum $\bigodot_{s\in P} X_s$ is dense in $F(X)$.

Suppose $P$ is left-cancellative.  Then for any $x\in X$ and $\oplus x_t\in F(X)$
we have $p(xx_s) = p(xx_t)$ if and only if $s = t$,
so there is an element $(y_s) \in\prod X_s$ such that
\[
y_s = \begin{cases} xx_t & \text{if $s = p(x)t$} \\
0 & \text{if $s\notin p(x)P$;} \end{cases}
\]
we write $(xx_t)$ for $(y_s)$.
Since
$\langle xx_s, xx_s \rangle_A \le \lVert x \rVert^2 \langle x_s, x_s \rangle_A$
for each $s\in P$,
the series $\sum \langle xx_s, xx_s \rangle_A$ is summable.
It is routine to check that
\[
l(x)(\oplus x_s) := \oplus xx_s
\qquad\text{for $\oplus x_s\in F(X)$}
\]
determines an adjointable operator $l(x)$ on $F(X)$;
indeed, the adjoint $l(x)^*$ is zero on any summand $X_s$
for which $s\notin p(x)P$, and on $X_{p(x)t} = \clsp X_{p(x)}X_t$
it is determined by the formula
$l(x)^*(yz) = \langle x, y\rangle_A\cdot z$ for $y\in X_{p(x)}$ and $z\in X_t$.
It follows that $l:X\to\Ll(F(X))$
is a Toeplitz representation, called the {\em Fock representation of $X$\/}.
The homomorphism $l_e:A\to\Ll(F(X))$ is simply the diagonal left action of $A$; that is,
$l_e(a) = \oplus \phi_s(a)$.  Since $\phi_e$ is just left multiplication on $X_e = A$,
it is isometric, and hence so is $l_e$;
by Remark~\ref{remark:toeplitz rep}, $l$ is isometric.

\begin{prop}\label{prop:toeplitz algebra}
Let $X$ be a product system over $P$ of Hilbert $A$--$A$ bimodules.
Then there is a $C^*$\ndash al\-ge\-bra $\Tt_X$,
called the {\em Toeplitz algebra\/} of $X$,
and a Toeplitz representation $i_X:X\to\Tt_X$, such that

\textup{(a)} for every Toeplitz representation $\psi$ of $X$,
there is a homomorphism $\psi_*$ of $\Tt_X$ such that
$\psi_*\circ i_X = \psi$; and

\textup{(b)} $\Tt_X$ is generated as a $C^*$\ndash al\-ge\-bra by $i_X(X)$.

\noindent The pair $(\Tt_X, i_X)$ is unique up to canonical isomorphism,
and $i_X$ is isometric.
\end{prop}

\begin{proof} It is straightforward to translate
the proof of \cite[Proposition~1.3]{fowrae2} to this setting.
\end{proof}

\begin{prop}\label{prop:CP algebra}
Let $X$ be a product system over $P$ of Hilbert $A$--$A$ bimodules.
Then there is a $C^*$\ndash al\-ge\-bra $\Oo_X$,
called the {\em Cuntz-Pimsner algebra\/} of $X$,
and a Toeplitz representation $j_X:X\to\Oo_X$ which is Cuntz-Pimsner covariant,
such that

\textup{(a)} for every Cuntz-Pimsner covariant Toeplitz representation $\psi$ of $X$,
there is a homomorphism $\psi_*$ of $\Oo_X$ such that
$\psi_*\circ j_X = \psi$; and

\textup{(b)} $\Oo_X$ is generated as a $C^*$\ndash al\-ge\-bra by $j_X(X)$.

\noindent The pair $(\Oo_X, j_X)$ is unique up to canonical isomorphism.
\end{prop}

\begin{remark}
Although the universal map $i_X:X\to\Tt_X$ is always isometric,
it is quite possible that $X$ might not admit any nontrivial
Toeplitz representations which are Cuntz-Pimsner covariant,
in which case $\Oo_X$ is trivial.
\end{remark}

\begin{proof}[Proof of Proposition~\ref{prop:CP algebra}]
With $(\Tt_X,i_X)$ as in Proposition~\ref{prop:toeplitz algebra},
let $\Ii$ be the ideal in $\Tt_X$ generated by
\[
\{ i_X(a) - i_X^{(s)}(\phi_s(a)) : s\in P, a \in \phi_s^{-1}(\Kk(X)) \}.
\]
Define $\Oo_X := \Tt_X/\Ii$ and $j_X := q\circ i_X$, where $q:\Tt_X \to \Oo_X$
is the canonical projection.  Obviously $j_X$ is a Toeplitz representation
which generates $\Oo_X$,
and it is Cuntz-Pimsner covariant because $j_X^{(s)} = q\circ i_X^{(s)}$.
If $\psi$ is another Cuntz-Pimsner covariant Toeplitz representation,
then the homomorphism $\psi_*$ of $\Tt_X$ satisfies
\[
\psi_*(i_X(a) - i_X^{(s)}(\phi_s(a))) = \psi(a) - \psi^{(s)}(\phi_s(a)) = 0
\]
whenever $\phi_s(a) \in \Kk(X_s)$,
and hence $\psi_*$ descends to the required homomorphism of $\Oo_X$
(also denoted $\psi_*$).
\end{proof}

\begin{prop}\label{prop:ps over N}
Let $X$ be a product system over $\NN$ of Hilbert $A$--$A$ bimodules.
Then $\Tt_X$ is canonically isomorphic to the Toeplitz algebra $\Tt_{X_1}$
of the Hilbert bimodule $X_1$.
If the left action on each fiber is isometric,
or if the left action on each fiber is by compact operators,
then $\Oo_X$ is canonically isomorphic to $\Oo_{X_1}$.
\end{prop}

\begin{proof}
Let $i_X:X\to\Tt_X$ be universal for Toeplitz representations of $X$,
and define $\mu := (i_X)_1 : X_1 \to \Tt_X$
and $\pi := (i_X)_0 : A = X_0  \to \Tt_X$.
Since $(\mu,\pi)$ is a  Toeplitz  representation of $X_1$,
we get a homomorphism $\mu\times\pi:\Tt_{X_1} \to \Tt_X$.

To construct the inverse of $\mu\times\pi$,
let $(i_{X_1}, i_A)$ be the universal Toeplitz representation of $X_1$ in $\Tt_{X_1}$,
and fix $n\ge 1$.  By \cite[Proposition~1.8(1)]{fowrae2},
there is a linear map $\psi_n:X_n\to\Tt_{X_1}$ which satisfies
\[
\psi_n(x_1\dotsm x_n) := i_{X_1}(x_1)\dotsm i_{X_1}(x_n)
\qquad\text{for all $x_1, \dots, x_n\in X_1$,}
\]
and then $(\psi_n, i_A)$ is a Toeplitz representation of $X_n$.
Defining $\psi_0 := i_A$ thus gives a Toeplitz representation $\psi:X\to\Tt_{X_1}$,
and it is routine to check that $\psi_*:\Tt_X\to\Tt_{X_1}$ is the inverse of $\mu\times\pi$.

Now let $i_X:X\to\Oo_X$ be universal for
Cuntz-Pimsner covariant Toeplitz representations of $X$.
As above, we get a homomorphism $\mu\times\pi:\Oo_{X_1}\to\Oo_X$.
To construct the inverse, we let $(i_{X_1},i_A):(X_1,A)\to \Oo_{X_1}$
be universal and define a Toeplitz representation $\psi:X\to\Oo_{X_1}$
as before; we need to check that $\psi$ is Cuntz-Pimsner covariant
under each of the hypotheses on the left action.
By definition $(\psi_1,\psi_0)$ is Cuntz-Pimsner covariant, so we use induction.
Assume that $(\psi_n,\psi_0)$ is Cuntz-Pimsner covariant for some $n\ge 1$,
and suppose $a\in A$ acts compactly on the left of $X_{n+1}$;
that is, $\phi(a)\otimes_A 1^n \in\Kk(X_{n+1})$.
Since the left action is isometric on each fiber,
by \cite[Lemma~7.2]{fmr} we have $\phi(a)\otimes_A 1^{n-1}\in\Kk(X_n)$;
hence $\psi^{(n)}(\phi(a)\otimes_A 1^{n-1}) = \psi_0(a)$.
But \cite[Lemma~7.5]{fmr} gives
$\psi^{(n+1)}(\phi(a)\otimes_A 1^n) = \psi^{(n)}(\phi(a)\otimes_A 1^{n-1})$,
so $(\psi_{n+1},\psi_0)$ is Cuntz-Pimsner covariant.

Now suppose that $A$ acts by compact operators on each fiber.
By representing $\Oo_{X_1}$ faithfully on a Hilbert space $\Hh$
we can assume that $\psi$ is a Toeplitz  representation of $X$ on $\Hh$.
Assuming again that $(\psi_n,\psi_0)$ is Cuntz-Pimsner covariant for some $n\ge 1$,
\cite[Lemma~1.6]{fmr} gives
$\overline{\psi_0(A)\Hh} \subseteq \clsp(\psi_n(X_n)\Hh)$.
Let $x\in X_n$, and express $x = y\cdot a$ with $y\in X_n$ and $a\in A$.
Since $(\psi_1,\psi_0)$ is Cuntz-Pimsner covariant and $\phi(a)\in\Kk(X_1)$, we have
\[
\psi_n(x) = \psi_n(y)\psi_0(a) = \psi_n(y)\psi^{(1)}(\phi(a)).
\]
Now $\phi(a)$ can be approximated by a finite sum $\sum \Theta_{x_i,y_i}$,
hence $\psi_n(x)$ can be approximated by a finite sum
$\psi_n(y)\psi(x_i)\psi(y_i)^* = \psi_{n+1}(yx_i)\psi(y_i)^*$.
Thus
\[
\overline{\psi_0(A)\Hh}
\subseteq \clsp(\psi_n(X_n)\Hh)
\subseteq\clsp(\psi_{n+1}(X_{n+1})\Hh),
\]
and $(\psi_{n+1}, \psi_0)$ is Cuntz-Pimsner covariant by \cite[Lemma~1.6]{fmr}.
\end{proof}

\begin{definition}\label{defn:nondegenerate}
Let $X$ be a product system over $P$ of Hilbert $A$--$A$ bimodules.
A Toeplitz representation $\psi:X\to B$ is {\em nondegenerate\/}
if the induced homomorphism $\psi_*:\Tt_X\to B$ is nondegenerate.
\end{definition}

\begin{lemma}\label{lemma:nondegenerate}
Suppose each fiber $X_s$ is essential as a left $A$-module.
Then a Toeplitz representation $\psi:X\to B$ is nondegenerate
if and only if the homomorphism $\psi_e:A\to B$ is nondegenerate.
\end{lemma}

\begin{proof}
Let $(a_i)$ be an approximate identity for $A = X_e$.
By \eqref{eq:approximate identity},
$i_X(a_i)$ is an approximate identity for $\Tt_X$,
and the result follows.
\end{proof}

\section{Crossed products twisted by multipliers}
\label{section:crossed products}

Our main examples of product systems come from $C^*$\ndash dynamical systems.
Suppose $\beta$ is an action of $P$ as endomorphisms of $A$
such that $\beta_e$ is the identity endomorphism.
We will assume that each $\beta_s$ is {\em extendible\/};
that is, that each $\beta_s$ extends to a strictly continuous endomorphism
$\overline{\beta_s}$ of $M(A)$.
For $P$ the positive cone of a totally ordered abelian group,
Adji has shown that extendibility is necessary to define a reasonable
crossed product $A\cross_\beta P$ \cite{adji}.

In this section we will consider crossed products which are twisted by
a multiplier $\omega$ of $P$; that is, by a function
$\omega:P\times P\to\TT$ which satisfies
$\omega(e,e) = 1$ and
\[
\omega(r,s)\omega(rs,t) = \omega(r,st)\omega(s,t)
\qquad\text{for all $r,s,t\in P$.}
\]
We call $(A,P,\beta,\omega)$ a {\em twisted semigroup dynamical system\/}.

\begin{definition}\label{defn:omega crossed product}
Let $B$ be a $C^*$\ndash al\-ge\-bra.
A function $V:P\to B$ is called an
{\em $\omega$\ndash representation\/} of $P$ if
\begin{equation}\label{eq:omegarep}
V_sV_t = \omega(s,t)V_{st}
\qquad\text{for all $s,t\in P$.}
\end{equation}
If  in addition each $V_s$ is an isometry  (resp.  partial isometry),
$V$ is called {\em isometric\/} (resp. {\em partial isometric\/})
$\omega$\ndash representation.
A {\em covariant representation\/} of $(A,P,\beta,\omega)$ on a Hilbert space $\Hh$
is a pair $(\pi,V)$ consisting of a nondegenerate representation $\pi:A\to B(\Hh)$
and an isometric  $\omega$\ndash representation $V:P\to B(\Hh)$ such that
\begin{equation}\label{eq:piVcovariance}
\pi(\beta_s(a)) = V_s\pi(a)V_s^*
\qquad\text{for all $s\in P$ and $a\in A$.}
\end{equation}
A {\em crossed product\/} for  $(A,P,\beta,\omega)$  is a triple
$(B,i_A,i_P)$ consisting of a $C^*$\ndash al\-ge\-bra $B$,
a nondegenerate homomorphism $i_A:A\to B$, and a map $i_P:P\to M(B)$
such that

\textup{(a)} if $\sigma$ is a nondegenerate representation of $B$,
then $(\sigma\circ i_A, \overline\sigma\circ i_P)$
is a covariant representation of $(A,P,\beta,\omega)$;

\textup{(b)} for every covariant representation $(\pi,V)$,
there is a representation $\pi\times V$ such that
$(\pi\times V)\circ i_A = \pi$ and $\overline{\pi\times V}\circ i_P = V$;
and

\textup{(c)} $B$ is generated as a $C^*$\ndash al\-ge\-bra by
$\{i_A(a)i_P(s): a\in A, s\in P\}$.
\end{definition}

After establishing the existence of a crossed product,
it is easily seen to be unique up to canonical isomorphism;
we denote the crossed product algebra $A\cross_{\beta,\omega} P$.

We will construct a product system $X = X(A,P,\beta,\omega)$
over the opposite semigroup $P^o$,
and show that its Cuntz-Pimsner algebra $\Oo_X$
is a crossed product for $(A,P,\beta,\omega)$.
Moreover, we will show that the Toeplitz  algebra of this product system
also has a crossed product structure:
it will be universal for pairs $(\pi,V)$ satisfying \eqref{eq:piVcovariance}
in which $\pi$ is a nondegenerate representation of $A$
and $V$ is a partial isometric $\omega$\ndash representation
such that
\begin{equation}\label{eq:Vpartialisometry}
\text{$V_s^*V_s\pi(a) = \pi(a)V_s^*V_s$ for all $s\in P$ and $a\in A$.}
\end{equation}
We call such a pair $(\pi,V)$ a {\em Toeplitz covariant representation\/} of $(A,P,\beta,\omega)$, 
and write $\Tt(A\cross_{\beta,\omega} P)$ for the
corresponding universal $C^*$\ndash al\-ge\-bra, called the {\em Toeplitz crossed product\/}
of $(A,P,\beta,\omega)$.

For each $s\in P$ let
\[
X_s := \{s\}\times \overline{\beta_s}(1)A,
\]
and give $X_s$ the structure of a Hilbert bimodule over $A$ via
\[
(s,x)\cdot a := (s,xa),
\quad
\langle (s,x),(s,y) \rangle_A := x^*y,
\]
and
\[
\phi_s(a)(s,x) := (s,\beta_s(a)x).
\]
Let $X = \bigsqcup_{s\in P} X_s$, let $p(s,x) := s$,
and define multiplication in $X$ by
\[
(s,x)(t,y) := (ts, \overline{\omega(t,s)}\beta_t(x)y )
\qquad\text{for $x\in\overline{\beta_s}(1)A$ and $y\in\overline{\beta_t}(1)A$.}
\]

\begin{lemma}\label{lemma:crossed product ps}
$X = X(A,P,\beta,\omega)$  is a product system  over the opposite semigroup $P^o$.
For each $s\in P$, the fiber $X_s$ is essential as a left $A$-module,
and the left action of $A$ on $X_s$ is by compact operators.
\end{lemma}

\begin{proof}  Let $(s,x)\in X_s$ and $(t,y) \in X_t$.
If $x = \overline{\beta_s}(1)a$ and $y = \overline{\beta_t}(1)b$, then
\[
\beta_t(x)y
= \beta_t(\overline{\beta_s}(1)a)\overline{\beta_t}(1)b
= \overline{\beta_t}(\overline{\beta_s}(1))\beta_t(a)b
= \overline{\beta_{ts}}(1)\beta_t(a)b,
\]
and hence the product $(s,x)(t,y)$ belongs to $X_{ts}$.
Letting $a$ vary over an approximate identity for $A$,
this product converges in norm to $\overline{\beta_{ts}}(1)b$,
so the set of products $(s,x)(t,y)$ has dense linear span in $X_{ts}$.
Hence to see that multiplication induces an isomorphism $X_s \otimes_A X_t \to X_{ts}$,
it suffices to check that it preserves the inner product of any pair of elementary
tensors:
\begin{align*}
\langle (s,x)\otimes_A (t,y),
&  (s,x')\otimes_A (t,y') \rangle_A \\
& = \langle (t,y), \phi_t(\langle (s,x), (s,x') \rangle_A)(t,y')\rangle_A \\
& = \langle (t,y), \phi_t(x^*x')(t,y') \rangle_A \\
& = \langle (t, y), (t, \beta_t(x^*x')y') \rangle_A \\
& = y^*\beta_t(x^*x')y' \\
& = \langle (ts, \overline{\omega(t,s)}\beta_t(x)y), 
            (ts, \overline{\omega(t,s)}\beta_t(x')y') \rangle_A \\
& = \langle (s,x)(t,y), (s,x')(t,y') \rangle_A.
\end{align*}
Multiplication is associative since
\begin{align*}
((s,x)(t,y))(r,z)
& = (ts, \overline{\omega(t,s)}\beta_t(x)y)(r,z) \\
& = (rts, \overline{\omega(r,ts)}\beta_r(\overline{\omega(t,s)}\beta_t(x)y)z) \\
& = (rts, \overline{\omega(rt,s)}\beta_{rt}(x)\overline{\omega(r,t)}\beta_r(y)z) \\
& = (s,x)(rt, \overline{\omega(r,t)}\beta_r(y)z) \\
& = (s,x)((t,y)(r,z)).
\end{align*}

If $(a_i)$ is an approximate identity for $A$, then for each $a\in A$ and $s\in P$
we have
$\lim \phi(a_i)(s,\overline{\beta_s}(1)a)
 = \lim (s,\beta_s(a_i)a) = (s,\overline{\beta_s}(1)a)$,
so each $X_s$ is essential.  If $a\in A$, then by writing $a = bc^*$
with $b,c\in A$, we see that
$\phi_s(a) = \Theta_{(s,\beta_s(b)),(s,\beta_s(c))} \in\Kk(X_s)$
is compact.
\end{proof}

\begin{lemma}\label{lemma:strict convergence}
Let $i_X:X\to\Tt_X$ be universal for Toeplitz representations of $X$,
and let $(a_i)$ be an approximate identity for $A$.
Then for each $s\in P$, $i_X(s,\beta_s(a_i))$ converges strictly in $M\Tt_X$.
\end{lemma}

\begin{proof}
Since each fiber $X_t$ is essential, any vector $\xi\in X_t$ can be written in the form
$\xi = \phi_t(a)\eta\cdot b$ with $a,b\in A$ and $\eta\in X_t$.
But then $i_X(\xi) = i_X(e,a)i_X(\eta)i_X(e,b)$,
and since elements  of the form $i_X(\xi)$ generate $\Tt_X$ as a $C^*$\ndash al\-ge\-bra,
the result  follows from the calculations
\begin{equation}\label{eq:strict1}
i_X(s,\beta_s(a_i))i_X(e,a)
= i_X(s, \beta_s(a_i)a) \to i_X(s, \overline{\beta_s}(1)a)
\end{equation}
and
\begin{equation}\label{eq:strict2}
i_X(e,a)i_X(s,\beta_s(a_i))
= i_X(s,\beta_s(aa_i)) \to i_X(s,\beta_s(a)).
\end{equation}
\end{proof}

Define $i_A:A\to\Tt_X$ by $i_A(a) := i_X(e,a)$,
and define $i_P:P\to M\Tt_X$ by
$i_P(s) := \lim i_X(s,\beta_s(a_i))^*$.

\begin{prop}\label{prop:crossed product}
$\Tt_X$ and $\Oo_X$ are canonically isomorphic
to $A\cross_{\beta,\omega} P$ and $\Tt(A\cross_{\beta,\omega} P)$, respectively.
More precisely, $(\Tt_X, i_A, i_P)$ is a Toeplitz crossed product for $(A,P,\beta,\omega)$,
and, with $q:\Tt_X\to\Oo_X$ the canonical projection,
$(\Oo_X, q\circ i_A, \overline q\circ i_P)$ is a crossed product for $(A,P,\beta,\omega)$.
\end{prop}

\begin{proof}
Taking $s = e$ in \eqref{eq:strict1} and \eqref{eq:strict2},
shows that $i_A(a_i)$ converges strictly to the identity in $M(\Tt_X)$;
that is, $i_A$ is nondegenerate.
For condition (a) of a Toeplitz crossed product,
suppose $\sigma$ is a nondegenerate representation of $\Tt_X$;
we must show that $(\pi,V) := (\sigma\circ i_A, \overline\sigma\circ i_P)$
is a Toeplitz covariant representation of $(A,P,\beta,\omega)$.
First note that $\pi$ is nondegenerate since $\sigma$ and $i_A$ are.
Equation \eqref{eq:strict2} gives
\begin{equation}\label{eq:AP}
i_A(a)i_P(s)^* = i_X(s,\beta_s(a))
\qquad\text{for  all $a\in A$ and $s\in P$,}
\end{equation}
so
\begin{align*}
i_P(s)i_A(a)i_P(s)^*
& = \lim i_X(s,\beta_s(a_i))^*i_X(s,\beta_s(a)) \\
& = \lim i_X(e,\beta_s(a_i^*a))
  = i_X(e, \beta_s(a)) = i_A(\beta_s(a)),
\end{align*}
and applying $\overline\sigma$ gives $V_s\pi(a)V_s^* = \pi(\beta_s(a))$.
In particular
\[
i_P(s)i_P(s)^*
= \lim i_P(s)i_A(a_i)i_P(s)^*
= \lim i_A(\beta_s(a_i)) = \overline{i_A}(\overline{\beta_s}(1))
\]
is a projection, so $i_P(s)$, and hence $V_s$, is a partial isometry.

To establish \eqref{eq:Vpartialisometry},
take any $a\in A$, write $a = bc^*$ with $b,c\in A$, and compute:
\begin{equation}\label{eq:domain commutes}
\begin{split}
i_A(bc^*)i_P(s)^*i_P(s)
& = \lim i_X(s,\beta_s(bc^*))i_X(s,\beta_s(a_i))^*
    \qquad \text{(by \eqref{eq:AP})} \\
& = \lim i_X(s,\beta_s(b))i_X(e,\beta_s(c^*))i_X(s,\beta_s(a_i))^*  \\
& = \lim i_X(s,\beta_s(b)) \bigl( i_X(s,\beta_s(a_i)i_X(e,\beta_s(c)) \bigr)^*  \\
& = \lim i_X(s,\beta_s(b)) i_X(s,\beta_s(a_i c))^*  \\
& =      i_X(s,\beta_s(b)) i_X(s,\beta_s(c))^*.
\end{split}
\end{equation}
Taking adjoints, interchanging $b$ and $c$,
and applying $\overline\sigma$ gives
$V_s^*V_s\pi(a) = \pi(a)V_s^*V_s$.

For every $s,t\in P$ we have
\begin{align*}
i_P(t)^*i_P(s)^*
& = (\lim_i j_X(t,\beta_t(a_i)))(\lim_j  j_X(s,\beta_s(a_j))) \\
& = \lim_i \lim_j j_X(st, \overline{\omega(s,t)}\beta_s(\beta_t(a_i))\beta_s(a_j)) \\
& = \lim_i j_X(st, \overline{\omega(s,t)}\beta_{st}(a_i)) \\
& = \overline{\omega(s,t)}i_P(st)^*;
\end{align*}
taking adjoints and applying $\overline\sigma$ gives
$V_sV_t = \omega(s,t)V_{st}$.
This completes the proof of condition (a).

For condition (b),
suppose $(\pi,V)$ is a Toeplitz covariant representation
of $(A,P,\beta,\omega)$ on a Hilbert space $\Hh$.
Define $\psi:X\to B(\Hh)$ by
\[
\psi(s,x) := V_s^*\pi(a).
\]
Since $\pi$ is nondegenerate and
$\pi(a) = \pi(\beta_e(a)) = V_e\pi(a)V_e^*$ for all $a\in A$,
$V_e$ is a coisometry.  Since $V_e^2 = \omega(e,e)V_e = V_e$,
we deduce that $V_e = 1$.
Thus
\begin{align*}
\psi(s,x)^*\psi(s,y)
& = \pi(x)^*V_sV_s^*\pi(y)
  = \pi(x^*\overline{\beta_s}(1)y) \\
& = V_e^*\pi(x^*y)
  \qquad\qquad\text{(since $y\in\overline{\beta_s}(1)A$ and $V_e = 1$)} \\
&  = \psi(e,x^*y)
  = \psi(\langle (s,x), (s,y) \rangle_A),
\end{align*}
and since we also have
\begin{align*}
\psi(s,x)\psi(t,y)
& = V_s^*\pi(x)V_t^*\pi(y) & \\
& = V_s^*\pi(x)V_t^*V_tV_t^*\pi(y) & \text{($V_t$ is a partial isometry)} \\
& = V_s^*V_t^*V_t\pi(x)V_t^*\pi(y) & \text{(by \eqref{eq:Vpartialisometry})} \\
& = (V_tV_s)^*\pi(\beta_t(x))\pi(y) & \\
& = \overline{\omega(t,s)}V_{ts}^*\pi(\beta_t(x)y) & \\
& = \overline{\omega(t,s)}\psi(ts,\beta_t(x)y) & \\
& = \psi((s,x)(t,y)),&
\end{align*}
$\psi$ is a Toeplitz representation of $X$.
Let $\pi\times V$ be the representation
$\psi_*:\Tt_X\to B(\Hh)$.  Then
\[
(\pi\times V)\circ i_A(a) = \psi_*\circ  i_X(e,a) = \psi(e,a) = V_e^*\pi(a) = \pi(a),
\]
and
\begin{align*}
\overline{\pi\times V}\circ i_P(s)\pi(a)
& = \psi_*(i_P(s)i_A(a)) & \\
& = \psi_*(i_X(s,\beta_s(a^*))^*) & \text{(by \eqref{eq:AP})} \\
& = \psi(s,\beta_s(a^*))^*
  = \pi(\beta_s(a))V_s & \\
& = V_s\pi(a)V_s^*V_s
  = V_sV_s^*V_s\pi(a)
  = V_s\pi(a); &
\end{align*}
since $\pi$ is nondegenerate, this implies that
$\overline{\pi\times V}\circ  i_P = V$,
as required.
For  condition (c), simply note that
$i_A(a)i_P(s) = i_X(s, \overline{\beta_s}(1)a^*)^*$,
and elements of this form generate $\Tt_X$.

We now show that $(\Oo_X, q\circ  i_A, \overline q\circ i_P)$
is a crossed  product for $(A,P,\beta,\omega)$.
Since $i_A$ and $q$ are nondegenerate, so is $q\circ i_A$.
If $\rho$ is a nondegenerate  representation of $\Oo_X$,
then $\sigma := \rho\circ  q$ is a nondegenerate representation
of $\Tt_X$.
Hence $(\pi,V) := (\rho\circ q \circ i_A, \overline\rho\circ\overline q \circ i_P)
= (\sigma \circ i_A, \overline\sigma\circ i_P)$ is a Toeplitz covariant  representation
of $(A,P,\beta,\omega)$.
To see that each $V_s$ is an isometry,
let $b,c\in A$.  Since $q\circ i_X$ is Cuntz-Pimsner covariant,
\eqref{eq:domain commutes} gives
\begin{align*}
q\circ i_A(bc^*)\overline q \circ i_P(s)^*\overline q \circ i_P(s)
& = q\circ i_X(s,\beta_s(b)) q\circ i_X(s,\beta_s(c))^* \\
& = (q\circ i_X)^{(s)}(\Theta_{(s,\beta_s(b)),(s,\beta_s(c))}) \\
& = (q\circ i_X)^{(s)}(\phi_s(bc^*)) \\
& = q\circ i_X(e,bc^*) = q\circ i_A(bc^*).
\end{align*}
Since $q\circ i_A$ is nondegenerate,
this implies that $\overline q\circ i_P(s)$, and hence $V_s$, is an isometry.
This gives condition (a) for a crossed product.
Condition (c) is obvious, so it remains only to verify (b).
Suppose $(\pi,V)$ is a covariant representation of $(A,P,\beta,\omega)$ on a Hilbert space $\Hh$,
and define $\psi(s,x) := V_s^*\pi(x)$ as before.  We have  already seen that
$\psi$ is a Toeplitz representation of $X$, and it is Cuntz-Pimsner covariant since,
for any $b,c\in A$,
\begin{align*}
\pi^{(s)}(\phi_s(bc^*))
& = \pi^{(s)}(\Theta_{(s,\beta_s(b)),(s,\beta_s(c))})
  = \psi(s,\beta_s(b))\psi(s,\beta_s(c))^* \\
& = V_s^*\pi(\beta_s(b))\pi(\beta_s(c^*))V_s
  = V_s^*V_s\pi(bc^*)V_s^*V_s = \pi(bc^*).
\end{align*}
Defining $\pi\times V := \psi_* :\Oo_X\to B(\Hh)$ gives condition (c).
\end{proof}

\section{Crossed products twisted by product systems}
\label{section:crossed products2}

Multipliers of $P$ correspond to product systems over $P$ of one-dimensional
Hilbert spaces: given a multiplier $\omega$,
one defines multiplication on $P\times\CC$ by $(s,z)(t,w)  := (st, \omega(s,t)zw)$.
In this section we consider twisted semigroup dynamical systems
in which the multiplier $\omega$ is replaced by a product system $X$ of Hilbert bimodules,
and we construct a crossed product which is ``twisted by $X$''.  
For this, we first need to see how semigroups of endomorphism
arise from Toeplitz representations of product systems. 

\begin{prop}\label{prop:alpha}
\textup{(1)} Let  $X$ be a Hilbert bimodule over $A$, and suppose $(\psi,\pi)$
is a Toeplitz representation of $X$ on a Hilbert space $\Hh$.
Then there is a unique endomorphism
$\alpha = \alpha^{\psi,\pi}$ of $\pi(A)'$ such that
\begin{equation}\label{eq:alpha}
\alpha(S)\psi(x) = \psi(x)S
\qquad\text{for all $S\in\pi(A)'$ and $x\in X$,}
\end{equation}
and such that $\alpha(1)$ vanishes on $(\psi(X)\Hh)^\perp$.

\textup{(2)} Let $X$ be a product system over $P$ of Hilbert $A$--$A$ bimodules
in which each fiber $X_s$ is essential as a left $A$-module.
Let $\psi$ be a nondegenerate Toeplitz representation of $X$ on a Hilbert space $\Hh$,
and let $\alpha^\psi_s$ be the endomorphism $\alpha^{\psi_s,\psi_e}$ above.
Then $\alpha^\psi:P\to\End(\psi_e(A)')$ is a semigroup homomorphism,
and $\alpha^\psi_e$ is the identity endomorphism.
\end{prop}

\begin{proof} (1) The uniqueness of $\alpha$ is obvious.
By \cite[Proposition~2.69]{rw}, there is a unital homomorphism
$S\in\pi(A)' \mapsto 1\otimes_A S\in\Ind\pi(\Ll(X))' \subseteq B(X\otimes_A\Hh)$
determined by
\[
1\otimes_A S(x\otimes_A h) = x\otimes_A Sh
\qquad\text{for $x\in X$ and $h\in\Hh$.}
\]
Let $U:X\otimes_A\Hh \to \Hh$ be the isometry which satisfies
$U(x\otimes_A h) = \psi(x)h$ (see the proof of \cite[Proposition~1.6(1)]{fowrae2}),
and define
\[
\alpha(S) := U(1\otimes_A S)U^*
\qquad\text{for all $S\in\pi(A)'$.}
\]
Then $\alpha$ is a homomorphism, and $\alpha(1) = UU^*$ vanishes on $(\psi(X)\Hh)^\perp$.
If $S\in\pi(A)'$ and $x\in X$, then for any $h\in\Hh$ we have
\[
\alpha(S)\psi(x)h
 = U(1\otimes_A S)(x\otimes_A h)
 = U(x\otimes_A Sh) = \psi(x)Sh,
\]
giving \eqref{eq:alpha}.

Since $\pi(a)\psi(x)h = \psi(\phi(a)x)h$,
the space $\clsp\{\psi(x)h: x\in X, h\in\Hh\}$ reduces $\pi$;
hence for any $S\in\pi(A)'$ and $a\in A$,
both $\alpha(S)\pi(a)$ and $\pi(a)\alpha(S)$ vanish on $(\psi(X)\Hh)^\perp$.
This and
\begin{align*}
\alpha(S)\pi(a)\psi(x)h
& = \alpha(S)\psi(\phi(a)x)h
  = \psi(\phi(a)x)Sh \\
& = \pi(a)\psi(x)Sh
  = \pi(a)\alpha(S)\psi(x)h
\end{align*}
show that $\alpha(\pi(A)') \subseteq \pi(A)'$.

\textup{(2)} 
Let  $s,t\in P$, and suppose $x\in X_s$ and $y\in X_t$.
Vectors of the form  $xy$ have dense linear span in $X_{st}$;
since $X_t$ is essential, this holds even when $s = e$ (see Remark~\ref{remark:essential}).
Since
\begin{align*}
\alpha^\psi_s(\alpha^\psi_t(S))\psi_{st}(xy)
& = \alpha^\psi_s(\alpha^\psi_t(S))\psi_s(x)\psi_t(y) \\
& = \psi_s(x)\alpha^\psi_t(S)\psi_t(y)
  = \psi_s(x)\psi_t(y)S
  = \psi_{st}(xy)S,
\end{align*}
we deduce that
\begin{equation}\label{eq:intertwine}
\alpha^\psi_s\circ\alpha^\psi_t(S)\psi_{st}(z) = \psi_{st}(z)S
\qquad\text{for all $S\in\psi_e(A)'$ and $z\in X_{st}$.}
\end{equation}
Once we show that
$\alpha^\psi_s\circ\alpha^\psi_t(1) = \alpha^\psi_{st}(1)$,
it follows from the uniqueness of $\alpha^\psi_{st}$
that $\alpha^\psi_s\circ\alpha^\psi_t = \alpha^\psi_{st}$.

From \eqref{eq:intertwine} we see that
$\alpha^\psi_s\circ\alpha^\psi_t(1) \ge \alpha^\psi_{st}(1)$.
Suppose that $\alpha^\psi_s\circ\alpha^\psi_t(1)f = f$;
we will show that $\alpha^\psi_{st}(1)f = f$, which
will complete the proof.
Since $f$ is in the range of $\alpha^\psi_s(1)$, it can be approximated
by a finite sum $\sum_i \psi_s(x_i)g_i$.  Then
\[
f \doteq \alpha^\psi_s\circ\alpha^\psi_t(1)\sum_i \psi_s(x_i)g_i
 = \sum_i \psi_s(x_i)\alpha^\psi_t(1)g_i.
\]
Now each $\alpha^\psi_t(1)g_i$ can be approximated by a finite
sum $\sum_j \psi_t(y_{ij})h_{ij}$, and then
\[
f
 \doteq \sum_{i,j} \psi_s(x_i)\psi_t(y_{ij})h_{ij}
 = \sum_{i,j} \psi_{st}(x_iy_{ij})h_{ij}.
\]
Thus $f$ can be approximated arbitrarily closely by a vector in the range
of $\alpha^\psi_{st}(1)$, and hence $\alpha^\psi_{st}(1)f = f$.

Since each $X_t$ is essential, the assumption that $\psi$ is nondegenerate
implies that $\psi_e$ is a nondegenerate representation of $A$.
Since $\alpha^\psi_e(S)\psi_e(a)h = \psi_e(a)Sh  = S\psi_e(a)h$
for all $a\in A$ and $h\in\Hh$, we have $\alpha^\psi_e(S) = S$
for all $S\in\psi_e(A)'$.
\end{proof}

Consider a twisted semigroup dynamical system $(C, P, \beta, X)$
in which $C$ is a separable $C^*$\ndash al\-ge\-bra,
$\beta:P\to\End C$ is an action of the semigroup $P$
as extendible endomorphisms of $C$,
and $X$ is a product system over $P$ of Hilbert $A$--$A$ bimodules.
We assume that $\beta_e$ is the identity endomorphism,
and that each fiber $X_s$ is essential as a left $A$-module.

\begin{definition}\label{defn:covariant pair}
 A {\em covariant representation\/} of $(C,P, \beta, X)$ on a Hilbert space $\Hh$
is a pair $(L, \psi)$ consisting of a nondegenerate representation
$L:C\to B(\Hh)$ and a nondegenerate Toeplitz representation $\psi:X\to B(\Hh)$
such that

\textup{(i)} $L(C) \subseteq \psi_e(A)'$, and

\textup{(ii)} $L\circ\beta_s = \alpha^\psi_s\circ L$ for $s\in P$.
\end{definition}

\begin{definition}\label{defn:crossed product}
A {\em crossed product\/} for $(C, P, \beta, X)$ is a triple
$(B, i_C, i_X)$ consisting of a $C^*$\ndash al\-ge\-bra $B$,
a nondegenerate homomorphism $i_C:C\to M(B)$, and a nondegenerate Toeplitz representation
$i_X:X\to M(B)$ such that

(a) there is a faithful nondegenerate representation $\sigma$ of $B$
such that $(\overline\sigma\circ i_C, \overline\sigma\circ i_X)$
is a covariant representation of $(C, P, \beta, X)$;

(b) for every covariant representation $(L,\psi)$ of $(C, P, \beta, X)$,
there is a representation $L\times\psi$ of $B$
such that $(\overline{L\times\psi})\circ i_C = L$
and $(\overline{L\times\psi})\circ i_X = \psi$;

(c) the $C^*$\ndash al\-ge\-bra $B$ is generated by
$\{i_C(c)i_X(x): c\in C,\ x\in X\}$.
\end{definition}

\begin{remark}
If each fiber $X_s$ has a finite basis $\{u_{s,1},\dots, u_{s,n(s)}\}$
(in the sense that $x = \sum_k u_{s,k}\cdot\langle u_{s,k}, x\rangle_A$
for every $x\in X_s$), it is not  hard to show that (a) is equivalent
to asking that $i_C(c)i_X(a) = i_X(a)i_C(c)$ for all $c\in C$ and $a\in A=X_e$,
and that
\[
i_C(\beta_s(c)) = \sum_k i_X(u_{s,k})i_C(c)i_X(u_{s,k})^*
\qquad\text{for all $s\in P$ and $c\in C$.}
\]
In this case,
$(\overline\sigma\circ i_C, \overline\sigma\circ i_X)$
will be a covariant representation of $(C, P, \beta, X)$
for {\em every\/}  nondegenerate representation $\sigma$ of $B$;
however, as demonstrated in \cite[Example~2.5]{fowrae}
for product systems of Hilbert spaces,
in general one cannot expect this to be the case.
\end{remark}

\begin{prop}\label{prop:existence of cp}
If $(C, P, \beta, X)$ has a covariant representation, then it has
a crossed product $(\cp, i_C, i_X)$
which is unique in the following sense:
if $(B, i_C', i_X')$ is another crossed product for $(C, P, \beta, X)$,
then there is an isomorphism $\theta: \cp \to B$
such that $\overline\theta\circ i_C = i_C'$ and $\theta\circ i_X = i_X'$.
\end{prop}

\begin{remark} When $X$ is the product system $P\times\CC$
with multiplication given by a multiplier $\omega$,
it is not hard to see that $\cp$ is precisely the crossed product
$C\cross_{\beta,\omega} P$ defined in the previous section.
If $C$ is unital and $A = \CC$, then $\cp$ is the crossed product
defined in \cite[Section~2]{fowrae}.
\end{remark}

\begin{proof}[Proof of Proposition~\ref{prop:existence of cp}]
If $S$ is a set of pairs $(L,\psi)$ consisting of maps $L:C\to B(\Hh_{L,\psi})$
and $\psi:X\to B(\Hh_{L,\psi})$, then $(\oplus L,\oplus\psi)$ is a covariant
representation of $(C,P,\beta,X)$ if and only if each $(L,\psi)$ is.
The main point here is that the value of $\alpha^{\oplus\psi}_s$
on an element of $(\oplus\psi)_e(A)'$ of the form $\oplus L(c)$
is $\oplus\alpha^\psi_s(L(c))$.

Suppose $(L,\psi)$ is a nondegenerate covariant representation
on a separable Hilbert space $\Hh$; that is, the $C^*$\ndash al\-ge\-bra
\[
\Uu := C^*(\{ L(c)\psi(x): c\in C,  x\in X\})
\]
acts nondegenerately on $\Hh$. 
We will identify the multiplier algebra of $\Uu$
with the concrete $C^*$\ndash al\-ge\-bra
\[
M(\Uu) = \{S\in B(\Hh): ST, TS\in\Uu\text{ for every $T\in\Uu$}\}.
\]
We claim that $L(C)\cup\psi(X)\subseteq M(\Uu)$.
For this, it suffices to check
that multiplying a generator $L(c)\psi(x)$ of $\Uu$ on either the left
or the right by an operator of the form $L(d)$, $\psi(y)$, or $\psi(y)^*$
yields another element of $\Uu$.
Certainly $L(d)L(c)\psi(x) = L(dc)\psi(x)\in\Uu$
and $L(c)\psi(x)\psi(y) = L(c)\psi(xy)\in\Uu$,
and since
\begin{equation}\label{eq:psiL}
\psi(y)L(c) = \alpha^\psi_{p(y)}(L(c))\psi(y) = L(\beta_{p(y)}(c))\psi(y),
\end{equation}
we also have $\psi(y)L(c)\psi(x) = L(\beta_{p(y)}(c))\psi(yx)\in\Uu$
and $L(c)\psi(x)L(d) = L(c\beta_{p(x)}(d))\psi(x)\in\Uu$.
Writing $c  = c_1^*c_2$ with $c_1, c_2\in C$ gives
\[
\psi(y)^*L(c)\psi(x) = (L(c_1)\psi(y))^*L(c_2)\psi(x) \in \Uu.
\]
Finally, to see that $L(c)\psi(x)\psi(y)^*\in\Uu$, we use  a trick from \cite{adji}.
Let $(c_i)$ be an approximate identity for $C$;
we claim that
\begin{equation}\label{eq:norm convergence}
L(c)\psi(x)L(c_i) \stackrel{\lVert\ \rVert}\longrightarrow L(c)\psi(x).
\end{equation}
Since $L$ is nondegenerate, $L(c)\psi(x)L(c_i)$
converges strongly to $L(c)\psi(x)$.
On the other hand, using \eqref{eq:psiL} we see that
$L(c)\psi(x)L(c_i)
= L(c\beta_{p(x)}(c_i))\psi(x)$
converges in {\em norm\/}
(to $L(c\overline{\beta_{p(x)}}(1))\psi(x)$),
and \eqref{eq:norm convergence} follows.
Hence
\[
L(c)\psi(x)L(c_i)\psi(y)^* \stackrel{\lVert\ \rVert}\longrightarrow L(c)\psi(x)\psi(y)^*,
\]
and since
\[
L(c)\psi(x)L(c_i)\psi(y)^*
= L(c)\psi(x)\psi(y)^*L(\beta_{p(y)}(c_i)) \in\Uu,
\]
we deduce that $L(c)\psi(x)\psi(y)^*\in\Uu$.

Since $M(\Uu) \subseteq \Uu''$, we have shown
in particular that the ranges of $L$ and $\psi$ are contained in $\Uu''$.
Consequently, any decomposition $1 = \sum Q_\lambda$ of the identity
as a sum of mutually orthogonal projections $Q_\lambda\in\Uu'$ gives
corresponding decompositions $L = \oplus Q_\lambda L$ and $\psi = \oplus Q_\lambda\psi$,
and by the first paragraph each pair $(Q_\lambda L,Q_\lambda\psi)$ is a covariant
representation of $(C,P,\beta,X)$.  By the usual Zorn's Lemma argument we
can choose these projections such that $\Uu$ acts cyclically on $Q_\lambda\Hh$;
since $C^*(\{Q_\lambda L(c) Q_\lambda\psi(x): c\in C,\ x\in X\}) = Q_\lambda\Uu$
acts cyclically on $Q_\lambda\Hh$, this shows that every covariant representation of
$(C,P,\beta,X)$ decomposes as a direct sum of cyclic representations.

Let $S$ be a set of cyclic covariant representations with
the property that every cyclic covariant representation of $(C,P,\beta,X)$
is unitarily equivalent to an element in $S$.
It can be shown that such a set $S$ exists
by fixing a Hilbert space $\Hh$ of sufficiently
large cardinality (depending on the cardinalities of $C$ and $X$)
and considering only representations on $\Hh$.
Note that $S$ is nonempty because the system has a covariant representation,
which has a cyclic summand.

Define $i_C := \bigoplus_{(L,\psi)\in S} L$
and $i_X := \bigoplus_{(L,\psi)\in S} \psi$,
and let $\cp$ be the $C^*$\ndash al\-ge\-bra generated by
$\{i_C(c)i_X(x): c\in C,\ x\in X\}$.
By the first paragraph, $(i_C,i_X)$ is a covariant representation of $(C,P,\beta,X)$,
and it is nondegenerate since each $(L,\psi)$ is.
We deduce that both $i_C$ and $i_X$ map into $M(\cp)$,
and that condition (a) for a crossed product is  satisfied by taking
$\sigma$ to be the identity representation.   Condition (c) is trivial,
and (b) holds because every covariant representation decomposes as a direct sum
of cyclic representations.
We need to show that $i_C:C\to M(\cp)$ and $i_X:X\to M(\cp)$ are nondegenerate. 
For this, let $c\in C$ and $x\in X$.
If $(a_i)$ is an approximate identity for $A = X_e$,
then by \eqref{eq:approximate identity} we have
$i_C(c)i_X(x)i_X(a_i)  = i_C(c)i_X(x\cdot a_i) \to i_C(c)i_X(x)$
and
$i_X(a_i)i_C(c)i_X(x)  = i_C(c)i_X(a_i)i_X(x) = i_C(c)i_X(\phi(a_i)x) \to i_C(c)i_X(x)$,
so $i_X$ is nondegenerate (Lemma~\ref{lemma:nondegenerate}).
If $(c_i)$ is an approximate identity for $C$,
then $i_C(c_i)i_C(c)i_X(x) = i_C(c_ic)i_X(x)  \to i_C(c)i_X(x)$,
and since $i_C$ is nondegenerate as a representation on Hilbert space,
\eqref{eq:norm convergence} gives $i_C(c)i_X(x)i_C(c_i) \to i_C(c)i_X(x)$.
Thus $i_C$ is nondegenerate.

For the uniqueness assertion, suppose  $(B, i_C', i_X')$ is another crossed product.
Condition (a) allows us to assume that  $(i_C,i_X)$ and $(i_C', i_X')$
are covariant representations of $(C,P,\beta,X)$ on Hilbert spaces $\Hh$ and $\Hh'$.
Condition (b) then gives a representation $i_C'\times i_X':\cp\to B(\Hh')$
whose image is contained in $B$ since  $i_C'\times i_X'(i_C(c)i_X(x))  = i_C'(c)i_X'(x)$.
Similarly one obtains a map $i_C\times i_X:B\to\cp$ which is obviously an inverse
for $i_C'\times i_X':\cp\to B$.
\end{proof}

If $P$ is a subsemigroup of a group $G$,
then there is a {\em dual coaction\/} of $G$ on $\cp$:

\begin{prop}\label{prop:coaction}
Suppose $(C,P,\beta,X)$ is a twisted system which has a covariant
representation. If $P$ is a subsemigroup of a group $G$, then there is an injective
coaction
\[
\delta:\cp \to (\cp) \otimes_{\min} C^*(G)
\]
such that
\[
\delta(i_C(c)i_X(x)) = i_C(c)i_X(x) \otimes i_G(p(x)).
\]
If $G$ is abelian, there is a strongly continuous action
$\widehat\beta$  of $\widehat G$ on $\cp$ such that
\[
\widehat\beta_\gamma(i_C(c)i_X(x)) = \gamma(p(x))i_C(c)i_X(x).
\]
\end{prop}

\begin{proof} We follow \cite[Proposition~2.7]{fowrae}.
Let $\sigma$ be a faithful nondegenerate representation $\sigma$ of $\cp$
such that $(L,\psi) := (\overline\sigma\circ i_C, \overline\sigma\circ i_X)$
is a covariant representation of $(C, P, \beta, X)$,
and let $U$ be a unitary representation of $G$ whose integrated form $\pi_U$
is faithful on $C^*(G)$.
We claim that $(L\otimes 1, \psi\otimes (U\circ p))$
is a covariant representation of $(C, P, \beta, X)$.
Most of the verifications are routine, so we check only that
\begin{equation}\label{eq:coaction}
L(\beta_s(c))\otimes 1 = \alpha^{\psi\otimes(U\circ p)}_s(L(c)\otimes 1)
\qquad\text{for all $s\in P$ and $c\in C$.}
\end{equation}
For this, we show that $L(\beta_s(c))\otimes 1$ satisfies the properties
which characterize $\alpha^{\psi\otimes(U\circ p)}_s(L(c) \otimes 1)$
(Proposition~\ref{prop:alpha}).
First, let $x\in X_s$;  we show that \eqref{eq:coaction}
holds on any vector in the range of
$(\psi\otimes(U\circ p))(x) = \psi_s(x)\otimes U_s$:
\begin{multline*}
(L(\beta_s(c))\otimes 1)(\psi_s(x)\otimes U_s)
  = \alpha^\psi_s(L(c))\psi_s(x) \otimes U_s
  = \psi_s(x)L(c)\otimes U_s \\
  = (\psi_s(x)\otimes U_s)(L(c)\otimes 1)
  = \alpha^{\psi\otimes(U\circ p)}_s(L(c)\otimes 1)(\psi_s(x)\otimes U_s).
\end{multline*}
Next, note that $\alpha^{\psi\otimes(U\circ p)}_s(1)$
is the projection onto
\begin{align*}
& \clsp\{(\psi\otimes(U\circ p))(x)\xi: x\in X_s, \xi\in\Hh_\sigma\otimes\Hh_U\} \\
& \qquad\qquad = \clsp\{\psi_s(x)h \otimes U_sk : x\in X_s, h\in\Hh_\sigma, k\in\Hh_U\} \\
& \qquad\qquad = \clsp\{\psi_s(x)h : x\in X_s, h\in\Hh_\sigma\} \otimes \Hh_U,
\end{align*}
which is precisely the range of $\alpha^\psi_s(1)\otimes 1$.
Since $L(\beta_s(c))\otimes 1 = \alpha^\psi_s(L(c)) \otimes 1$
vanishes on the range of $1 - \alpha^\psi_s(1)\otimes 1$,
\eqref{eq:coaction} follows from the uniqueness assertion
of Proposition~\ref{prop:alpha}.

Since $(L\otimes 1, \psi\otimes (U\circ p))$ is covariant,
there is a representation $\rho$ of $\cp$ such that
\begin{align*}
\rho(i_C(c)i_X(x))
& = (L(c)\otimes 1)(\psi(x)\otimes U_{p(x)}) \\
& = (\sigma\otimes\pi_U)(i_C(c)i_X(x) \otimes i_G(p(x))).
\end{align*}
Since $\sigma$ and $\pi_U$ are faithful, $\sigma\otimes\pi_U$
is faithful on $(\cp) \otimes_{\min} C^*(G)$,
and we can define $\delta := (\sigma\otimes\pi_U)^{-1}\circ\rho$.

By checking on generators it is easy to see that $\delta$
satisfies the coaction identity
$(\id\otimes\delta_G) \circ \delta = (\delta\otimes\id) \circ \delta$,
and $\delta$ is injective since $\sigma = (\sigma\otimes\epsilon)\circ\delta$,
with $\epsilon$ the augmentation representation of $C^*(G)$
(i.e., $\epsilon(i_G(s)) = 1$ for all $s\in G$).
When $G$ is abelian, $\widehat\beta$ is the action
canonically associated  with $\delta$.
\end{proof}

\section{Nica covariance}
\label{section:covariance}

Now suppose $P$ is a subsemigroup of a group $G$ such that $P\cap P^{-1} = \{e\}$.
Then $s\le t$ iff $s^{-1}t\in P$ defines a partial order on $G$ which is left-invariant:
for any $r,s,t\in P$ we have $s\le t$ iff $rs \le rt$.
Following Nica \cite{nica}, we say that $(G,P)$ is a {\em quasi-lattice ordered group\/}
if every finite subset of $G$ which has an upper bound in $P$
has a least upper bound in $P$.
When $s,t\in P$ have a common upper bound, we denote their
least upper bound by $s\vee t$; when $s$ and $t$ have no common
upper bound we write $s\vee t = \infty$.
For a finite subset $C = \{t_1,\dots, t_n\}$ of $P$, we write
$\sigma C$ for $t_1 \vee \dots \vee t_n$.

\begin{definition}\label{defn:psicovariant}
Suppose $(G,P)$ is a quasi-lattice ordered group
and $X$ is a product system over $P$ of essential Hilbert $A$--$A$ bimodules.
We call a Toeplitz representation $\psi:X\to B(\Hh)$ {\em Nica covariant\/} if
\[
\alpha^\psi_s(1)\alpha^\psi_t(1)
= \begin{cases} \alpha^\psi_{s\vee t}(1) & \text{if $s\vee t < \infty$} \\
                0 & \text{otherwise.}
\end{cases}
\] 
\end{definition}

\begin{remark} If $(G,P)$ is totally ordered,
then every Toeplitz representation of $X$ is Nica covariant.
\end{remark}

\begin{lemma}\label{lemma:fock covariant}
Let $l:X\to\Ll(F(X))$ be the Fock representation, and suppose $\pi$
is a representation of $A$ on a Hilbert space $\Hh$.
Then
\[
\Psi := F(X)\dashind_A^{\Ll(X)}\pi\circ l
\]
is a Nica-covariant Toeplitz representation of $X$.
If $\pi$ is faithful, then $\Psi$ is isometric.
\end{lemma}

\begin{proof} Since $l$ is a Toeplitz representation, so is $\Psi$.
Let $s\in P$.  The range of $\alpha^\Psi_s(1)$ is
\begin{multline*}
\clsp\{\Psi(x)\xi: x\in X_s, \xi\in F(X)\otimes_A\Hh\} \\
 = \clsp\{l(x)y\otimes_A h:  x\in X_s, y\in F(X), h\in\Hh\}
 = \bigoplus_{s \le r} X_r \otimes_A \Hh.
\end{multline*}
Hence for any $s,t\in P$, the range of $\alpha^\Psi_s(1)\alpha^\Psi_t(1)$ is
\[
\Bigl( \bigoplus_{s \le r} X_r \otimes_A \Hh \Bigr)
\cap
\Bigl( \bigoplus_{t \le r} X_r \otimes_A \Hh \Bigr),
\]
which is $\bigoplus_{s\vee t \le r} X_r \otimes_A \Hh = \ran\alpha^\Psi_{s\vee t}(1)$
if $s\vee t < \infty$, and is zero otherwise.
If $\pi$ is faithful then so is $F(X)\dashind_A^{\Ll(X)}\pi$;
since $l$ is isometric, this implies that $\Psi$ is isometric.
\end{proof}

\begin{prop}\label{prop:CP Nica}
Let $(G,P)$ be a quasi-lattice ordered group
such that every $s,t\in P$ have a common upper bound.
Let $X$ be a product system over $P$ of essential Hilbert $A$--$A$ bimodules
such that the left action of $A$ on each fiber $X_s$ is by compact operators.
Then every Toeplitz representation $\psi:X\to B(\Hh)$
which is Cuntz-Pimsner covariant is also Nica covariant.
\end{prop}

\begin{proof} Fix $s\in P$.
Since $(\psi_s,\psi_e)$ is Cuntz-Pimsner covariant and $\phi_s(A)\subseteq\Kk(X_s)$,
\cite[Lemma~1.6]{fmr} gives $\psi_e(A)\Hh\subseteq\clsp\psi_s(X_s)\Hh$.
But $X_s$ is essential, so the reverse inclusion holds as well,
and since $\clsp\psi_s(X_s)\Hh$ is precisely the range of $\alpha^\psi_s(1)$,
we deduce that $\alpha^\psi_s(1)$ is constant in $s$.
Since $s\vee t<\infty$ for all $s,t\in P$, this implies that $\psi$
is  Nica covariant.
\end{proof}

There are product systems for which Nica covariance is not a
$C^*$\ndash al\-ge\-braic condition;
that is, if $\psi:X\to B(\Hh)$ is Nica covariant
and $\sigma:C^*(\psi(X))\to B(\Kk)$ is a homomorphism,
the composition $\sigma\circ\psi$ need not be Nica covariant
\cite[Example~1.3]{fowler}.
We pause a moment to show how to adapt the methods
of \cite{fowler} to avoid this pathology.
The following Lemma collects some results we shall need
for both this and the sequel.

\begin{lemma}\label{lemma:rho}
Suppose $X$ is a product system over $P$ of essential Hilbert $A$--$A$ bimodules,
$\psi:X\to B(\Hh)$ is a Toeplitz representation, and $s\in P$.

\textup{(1)}
There is a strict--strong continuous representation
$\rho^\psi_s:\Ll(X_s)\to B(\Hh)$
such that
\[
\rho^\psi_s(S)\psi_s(x)h  = \psi_s(Sx)h
\qquad\text{for all $S\in\Ll(X_s)$, $x\in X_s$, and $h\in\Hh$,}
\]
and such that $\rho^\psi_s(S)$ vanishes on $(\psi_s(X_s)\Hh)^\perp$.
Moreover, $\rho^\psi_s(S) = \psi^{(s)}(S)$ for every $S\in\Kk(X_s)$.

\textup{(2)} $\rho^\psi_s(1) = \alpha^\psi_s(1)$

\textup{(3)} If $a\in A$ satisfies $\phi_s(a)\in\Kk(X_s)$, then
\begin{equation}\label{eq:absorb}
\psi_e(a)\rho^\psi_s(1) = \psi^{(s)}(\phi_s(a)) = \rho^\psi_s(1)\psi_e(a).
\end{equation}

\textup{(4)} If $Q\in\psi_e(A)'$,
then $\alpha^\psi_s(Q)\in\rho^\psi_s(\Ll(X_s))'$.
Further, if $Q$ is a projection such that $\psi_e$ acts faithfully on $Q\Hh$,
then $\rho^\psi_s$ acts faithfully on $\alpha^\psi_s(Q)\Hh$.

\textup{(5)} For all $S\in\Ll(X_s)$ and $t\in P$ we have
\[
\rho^\psi_{st}(S\otimes_A 1)
= \rho^\psi_s(S)\rho^\psi_{st}(1)
= \rho^\psi_{st}(1)\rho^\psi_s(S),
\]
where $S\otimes_A 1(xy) := (Sx)y$ for all $x\in X_s$ and $y\in X_t$.

\textup{(6)} If $t\in P$ and $z,w\in X_s$, then
$\rho^\psi_{st}(\Theta_{z,w} \otimes_A 1) = \psi(z)\alpha^\psi_t(1)\psi(w)^*$.
\end{lemma}

\begin{proof} (1) See \cite[Proposition~1.6(1)]{fowrae2}.
For the continuity assertion,
suppose  $S_\lambda\to S$ strictly in $\Ll(X_s) = M\Kk(X_s)$,
$x\in X_s$, and $h\in\Hh$.
There exists $K\in\Kk(X_s)$ and $y\in X_s$ such that $x = Ky$,
and then
\[
\rho^\psi_s(S_\lambda)\psi_s(x)h
= \rho^\psi_s(S_\lambda K)\psi_s(y)h
\to \rho^\psi_s(SK)\psi_s(y)h
= \rho^\psi_s(S)\psi_s(x)h.
\]

(2) Both $\rho^\psi_s(1)$ and $\alpha^\psi_s(1)$ are the projection onto
$\clsp\{\psi_s(X_s)\Hh\}$.

(3) If $x\in X_s$ and $h\in\Hh$, then
\[
\psi_e(a)\rho^\psi_s(1)\psi_s(z)h
= \psi_e(a)\psi_s(z)h
= \psi_s(\phi_s(a)z)h
= \psi^{(s)}(\phi_s(a))\psi_s(z)h;
\]
since both sides of \eqref{eq:absorb} are supported on $\clsp\psi_s(X_s)\Hh$,
this implies that $\psi_e(a)\rho^\psi_s(1) = \psi^{(s)}(\phi_s(a))$.
By (2), $\rho^\psi_s(1)$ commutes with $\psi_e(a)$,
giving the other half of \eqref{eq:absorb}.

(4) When $Q$ is a projection, $\alpha^\psi_s(Q)$ is the projection
onto $\clsp\psi_s(X_s)Q\Hh$,
and the result follows from \cite[Proposition~1.6(2)]{fowrae2}.

(5) See \cite[Proposition~1.8(2)]{fowrae2}.

(6) $\rho^\psi_{st}(\Theta_{z,w} \otimes_A 1)
     = \rho^\psi_{st}(1)\rho^\psi_s(\Theta_{z,w})
     = \alpha^\psi_{st}(1)\psi(z)\psi(w)^*
     = \psi(z)\alpha^\psi_t(1)\psi(w)^*$.
\end{proof}

\begin{prop}\label{prop:characterize covariance}
Suppose $(G,P)$ is a quasi-lattice ordered group
and $X$ is a product system over $P$ of essential Hilbert $A$--$A$ bimodules.
A Toeplitz representation $\psi:X\to B(\Hh)$ is Nica covariant if and only if
\begin{equation}\label{eq:characterize covariance}
\rho^\psi_s(S)\rho^\psi_t(T)
= \begin{cases} \rho^\psi_{s\vee t}((S\otimes_A 1)(T\otimes_A 1))
    & \text{if $s\vee t < \infty$} \\
  0 & \text{otherwise.}
\end{cases}
\end{equation}
\end{prop}

\begin{proof} \cite[Proposition~1.4]{fowler}
If $\psi$ is Nica covariant, then
\begin{align*}
\rho^\psi_s(S)\rho^\psi_t(T)
& = \rho^\psi_s(S)\rho^\psi_s(1)\rho^\psi_t(1)\rho^\psi_t(T) \\
& = \rho^\psi_s(S)\rho^\psi_{s\vee t}(1)\rho^\psi_t(T)
  = \rho^\psi_{s\vee t}((S\otimes_A 1)(T\otimes_A 1)),
\end{align*}
where the last equality uses Lemma~\ref{lemma:rho}(5).
Conversely, suppose \eqref{eq:characterize covariance} holds for all compact $S$ and $T$.
If $S\to 1$ strictly, then
\[
\rho^\psi_{s\vee t}((S\otimes_A 1))
=   \rho^\psi_s(S)\rho^\psi_{s\vee t}(1)
\to \rho^\psi_s(1)\rho^\psi_{s\vee t}(1)
= \rho^\psi_{s\vee t}(1),
\]
where the convergence is in the strong operator topology.
Hence
\[
\rho^\psi_s(1)\rho^\psi_t(T)
= \begin{cases} \rho^\psi_{s\vee t}(T\otimes_A 1)
    & \text{if $s\vee t < \infty$} \\
  0 & \text{otherwise}
\end{cases}
\]
for every $T\in\Kk(X_t)$.
Letting $T\to 1$ strictly shows that $\psi$ is Nica covariant.
\end{proof}

When each product $(S\otimes_A 1)(T\otimes_A 1)$ is compact,
the previous Proposition allows us to give a $C^*$\ndash al\-ge\-braic
characterization of Nica covariance:

\begin{definition}\label{defn:CA}
Suppose $(G,P)$ is a quasi-lattice ordered group and
$X$ is a product system over $P$ of essential Hilbert $A$--$A$ bimodules.
We say that $X$ is {\em compactly aligned\/}
if whenever $s,t\in P$ have a common upper
bound and $S$ and $T$ are compact operators on $X_s$ and $X_t$, respectively,
$(S\otimes_A 1)(T\otimes_A 1)$
is a compact operator on $X_{s \vee t}$.
If $X$ is compactly aligned and $\psi$ is a Toeplitz representation of $X$ in
a $C^*$\ndash al\-ge\-bra $B$, we say that $\psi$ is {\em Nica covariant\/} if
\[
\psi^{(s)}(S)\psi^{(t)}(T)
= \begin{cases} \psi^{(s\vee t)}((S\otimes_A 1)(T\otimes_A 1))
    & \text{if $s\vee t < \infty$} \\
  0 & \text{otherwise.}
\end{cases}
\]
whenever $s,t\in P$, $S\in\Kk(X_s)$ and $T\in\Kk(X_t)$.
\end{definition}

\begin{prop}\label{prop:CA conditions}
If $(G,P)$ is a total order, or if the left action of $A$ on each fiber $X_s$
is by compact operators, then $X$ is compactly aligned.
\end{prop}

\begin{proof} Suppose $s,t\in P$, $s\vee t < \infty$,
$S\in\Kk(X_s)$, and $T\in\Kk(X_t)$.
If $(G,P)$ is a total order then either $S\otimes_A 1 = S$
or $T\otimes_A 1 = T$; either way $(S\otimes_A 1)(T\otimes_A 1)$
is compact.
If the left action of $A$ on each fiber $X_s$
is by compact operators, then by \cite[Corollary~3.7]{pimsner},
both $S\otimes_A 1$
and $T\otimes_A 1$ are compact.
\end{proof}

\begin{prop}\label{prop:preserve covariance}
Suppose $X$ is compactly aligned.
Let $B$ and $C$ be $C^*$\ndash al\-ge\-bras,
let $\psi:X\to B$ be a Nica-covariant Toeplitz representation,
and let $\sigma:B\to C$ be a homomorphism.
Then $\sigma\circ\psi$ is Nica covariant.
\end{prop}

\begin{proof} By checking on an operator  $\Theta_{x,y}\in\Kk(X_s)$,
one verifies that $(\sigma\circ\psi)^{(s)} = \sigma\circ\psi^{(s)}$,
and the result follows easily from this.
\end{proof}

\begin{prop}\label{prop:yz}
Suppose $X$ is a compactly-aligned product system,
$\psi$ is a Nica-covariant Toeplitz representation of $X$,
$s,t\in P$, $y\in X_s$, and $z\in X_t$.
If $s\vee t = \infty$, then $\psi(y)^*\psi(z) = 0$;
otherwise
\[
\psi(y)^*\psi(z)
\in\clsp\{\psi(f)\psi(g)^*: f\in X_{s^{-1}(s\vee t)}, g\in X_{t^{-1}(s\vee t)}\}.
\]
\end{prop}

\begin{proof} Express $y = Sy'$ with $S\in\Kk(X_s)$ and $y'\in X_s$;
similarly, express $z = Tz'$ with $T\in\Kk(X_t)$ and $z'\in X_t$.
Since $\psi$ is Nica covariant,
\[
\psi(y)^*\psi(z) = \psi(y')^*\rho^\psi_s(S^*)\rho^\psi_t(T)\psi(z')
\]
is zero if $s\vee t = \infty$, and otherwise
\[
\psi(y)^*\psi(z) = \psi(y')^*\rho^\psi_{s\vee t}(K)\psi(z'),
\]
where $K = (S^*\otimes_A 1)(T\otimes_A 1) \in\Kk(X_{s\vee t})$.
Since  $K$ is compact it can be approximated in norm by a finite sum of operators
$\Theta_{u,v}$ with $u,v\in X_{s\vee t}$,
and hence $\rho^\psi_{s\vee t}(K)$ can be approximated by finite sums
of the form $\psi(u)\psi(v)^*$.
But any such $u$ can be approximated by finite sums
of products $u_1f'$ with $u_1\in X_s$ and $f'\in X_{s^{-1}(s\vee t)}$;
similarly, any such $v$ can be approximated by finite sums
of products $v_1g'$ with $v_1\in X_t$ and $g'\in X_{t^{-1}(s\vee t)}$.
Hence $\psi(y')^*\rho^\psi_{s\vee t}(K)\psi(z')$ can be approximated in norm
by finite sums of operators of the form
\[
\psi(y')^*\psi(u_1)\psi(f')\psi(g')^*\psi(v_1)^*\psi(z')
= \psi(\langle y',u_1 \rangle_A f')\psi(\langle z',v_1 \rangle_A g')^*.
\]
\end{proof}

The following Lemma is useful
when working with Nica-covariant Toeplitz representations.

\begin{lemma}\label{lemma:handy}
Suppose $(G,P)$ is a quasi-lattice ordered group,
$X$ is a product system over $P$ of essential Hilbert $A$--$A$ bimodules,
$\psi$ is a Toeplitz representation of $X$ on $\Hh$,
$x\in X$, and $s\in P$.

\textup{(1)} If $p(x) \le s$, then
$\alpha^\psi_s(S)\psi(x) = \psi(x)\alpha^\psi_{p(x)^{-1}s}(S)$
for all $S\in\psi_e(A)'$.

\textup{(2)} If $\psi$ is Nica covariant, then
\[
\alpha^\psi_s(1)\psi(x) =
\begin{cases}
  \psi(x)\alpha^\psi_{p(x)^{-1}(p(x) \vee s)}(1)
  & \text{if $p(x) \vee s < \infty$,} \\
0 & \text{otherwise.}
\end{cases}
\]
\end{lemma}

\begin{proof} The proof is formally identical to that of \cite[Lemma~3.6]{fowrae}.
\end{proof}

\section{The system $(B_P,P,\tau,X)$}
\label{section:system}

For each $t\in P$, let $1_t\in\ell^\infty(P)$
be the characteristic function of $tP$.
Since the product $1_s1_t$ is either $1_{s\vee t}$ or $0$,
$B_P := \clsp\{1_t: t\in P\}$
is a $C^*$\ndash subalgebra of $\ell^\infty(P)$.
Left translation on $\ell^\infty(P)$ restricts to an action
$\tau$ of $P$ on $B_P$, determined by
$\tau_s(1_t) = 1_{st}$ for $s,t\in P$.

\begin{prop}\label{prop:Lpsi}
Suppose $(G,P)$ is a quasi-lattice ordered group and
$X$ is a product system over $P$ of essential Hilbert $A$--$A$ bimodules.

\textup{(1)} If $(L, \psi)$ is a covariant representation of
$(B_P, P, \tau, X)$, then $\psi$ is a nondegenerate Nica-covariant Toeplitz representation of $X$ and
$L(1_s) = \alpha^\psi_s(1)$.

\textup{(2)} If $\psi$ is a
nondegenerate Nica-covariant Toeplitz representation of $X$ on a Hilbert space $\Hh$,
then there is a representation $L^\psi:B_P\to B(\Hh)$ such that
$L^\psi(1_s) = \alpha^\psi_s(1)$; moreover, $(L^\psi, \psi)$ is 
then a covariant representation of $(B_P, P, \tau, X)$.
\end{prop}

\begin{proof}
The proof is formally identical to that of \cite[Proposition~4.1]{fowrae},
except that in (2) one must  also note that $L^\psi(B_P) \subseteq\psi_e(A)'$
since $L^\psi(1_s) = \alpha^\psi_s(1) \in \psi_e(A)'$
and $\{1_s: s\in P\}$ generates $B_P$.
\end{proof}

\begin{cor}\label{cor:Lpsi}
The system $(B_P, P, \tau, X)$ has a covariant representation.
\end{cor}

\begin{proof}
Let $\pi$ be a nondegenerate representation of $A$ on a Hilbert space $\Hh$,
and let $l:X\to\Ll(F(X))$ be the Fock representation of $X$.
By Lemma~\ref{lemma:fock covariant},
$\Psi := F(X)\dashind_A^{\Ll(X)}\pi\circ l$
is a Nica-covariant Toeplitz representation of $X$.
Since $\pi$ is nondegenerate, so is $F(X)\dashind_A^{\Ll(X)}\pi$;
since $l$ is nondegenerate, $\Psi$ is as well.
The previous Proposition thus gives a covariant representation
$(L^\Psi,\Psi)$ of $(B_P,P,\tau,X)$.
\end{proof}

Let $i_X$ and $i_{B_P}$ be the canonical maps of $X$ and $B_P$ into $\bpp$.
Since $B_P$ is unital, $i_X(x) = i_{B_P}(1)i_X(x) \in\bpp$ for each $x\in X$.
We write $\Tt_{\cv}(X)$ for the $C^*$\ndash sub\-algebra of $\bpp$ generated
by $i_X(X)$; the following Theorem justifies this notation.

\begin{theorem}\label{theorem:subalgebra}
$(\Tt_{\cv}(X),i_X)$ is universal for Nica-covariant Toeplitz representations of $X$,
in the sense that:

\textup{(a)} there is a faithful representation $\theta$ of $\Tt_{\cv}(X)$
on Hilbert space such that $\theta\circ i_X$ is a Nica-covariant
Toeplitz representation of $X$, and

\textup{(b)} for every Nica-covariant Toeplitz representation $\psi$ of $X$,
there is a representation $\psi_*$ of $\Tt_{\cv}(X)$ such that $\psi=\psi_*\circ i_X$.

\noindent Up to canonical isomorphism, $(\Tt_{\cv}(X),i_X)$ is the unique
pair with this property.
If $X$ is compactly aligned, then $i_X$ is Nica covariant,
\begin{equation}\label{eq:span toeplitz}
\Tt_{\cv}(X) = \clsp\{i_X(x)i_X(y)^*: x,y\in X\},
\end{equation}
and
\begin{equation}\label{eq:span bpp}
\bpp = \clsp\{i_X(x)i_{B_P}(1_s)i_X(y)^*: x,y\in X, s\in P\}.
\end{equation}
\noindent If the left  action of $A$ on each fiber $X_s$ is by compact operators,
then $\Tt_{\cv}(X)$ is all of $\bpp$;
if in addition every $s,t\in P$ have a common upper bound,
then the Cuntz-Pimsner algebra $\Oo_X$ is a quotient of $\Tt_{\cv}(X)$.
\end{theorem}

\begin{proof}[Proof of Theorem~\ref{theorem:subalgebra}]
Let $\sigma$ be a faithful representation of $\bpp$ on a Hilbert space $\Hh$
such that $(\overline\sigma \circ i_{B_P}, \sigma\circ i_X)$
is a covariant representation of $(B_P, P, \tau, X)$.
By Proposition~\ref{prop:Lpsi}(1), $\sigma\circ i_X$
is a Nica-covariant Toeplitz representation of $X$,
so  we can take  $\theta$ to be the restriction  of $\sigma$ to $\Tt_{\cv}(X)$.
Suppose $\psi$ is a (nondegenerate) Nica-covariant Toeplitz representation of $X$.
Proposition~\ref{prop:Lpsi}(2) gives us a covariant representation
$(L^\psi,\psi)$ of $(B_P,P,\tau,X)$,
and hence a representation $L^\psi\times\psi$ of $\bpp$ such that
$(L^\psi\times\psi)\circ i_X=\psi$.  Restricting
$L^\psi\times\psi$ to $\Tt_{\cv}(X)$ gives the required representation $\psi_*$.
Uniqueness of $(\Tt_{\cv}(X),i_X)$ follows by the usual argument.

Suppose $X$ is compactly aligned.  Since $i_X$ is the composition of the Nica-covariant
Toeplitz representation $\sigma\circ i_X$ and the homomorphism $\sigma^{-1}$
(restricted to $\sigma(\Tt_{\cv}(X))$),
$i_X$ is Nica covariant by Proposition~\ref{prop:preserve covariance}.
Let $w\in X$, and express $w = z\cdot a$
for some $z\in X$ and $a\in A$.
Then $i_X(w) = i_X(z)i_X(a^*)^*$, so
$\Aa := \clsp\{i_X(x)i_X(y)^*: x,y\in X\}$
contains $i_X(X)$.
Obviously $\Aa$ is a closed self-adjoint subspace of $\Tt_{\cv}(X)$,
and since $X$ is compactly aligned, Proposition~\ref{prop:yz} shows that
$\Aa$ is closed under multiplication.  This gives \eqref{eq:span toeplitz}.

Now let $\Bb := \clsp\{i_X(x)i_{B_P}(1_s)i_X(y)^*: x,y\in X, s\in P\}$.
Using Lemma~\ref{lemma:handy} with $\psi := \sigma\circ i_X$,
and then applying $\sigma^{-1}$, gives
\begin{equation}\label{eq:handy1}
i_X(x)i_{B_P}(1_s) = i_{B_P}(1_{p(x)s})i_X(x)
\end{equation}
and
\begin{equation}\label{eq:handy2}
i_{B_P}(1_s)i_X(x) = 
\begin{cases} i_X(x)i_{B_P}(1_{p(x)^{-1}(p(x)\vee s)}) & \text{if $p(x)\vee s < \infty$} \\
  0 & \text{otherwise.} \end{cases}
\end{equation}
Equation \eqref{eq:handy1} shows that
\[
i_X(x)i_{B_P}(1_s)i_X(y)^*
= i_{B_P}(1_{p(x)s})i_X(x)(i_{B_P}(1_{p(y)s})i_X(y))^*
\in\bpp,
\]
so $\Bb\subseteq\bpp$.
Since $B_P$ is generated by $\{1_s: s\in P\}$,
elements of the form $i_{B_P}(1_s)i_X(w)$ generate $\bpp$ as a $C^*$\ndash al\-ge\-bra;
with $w = z\cdot a$ as above, \eqref{eq:handy2} shows that
\begin{align*}
i_{B_P}(1_s)i_X(w)
& = i_{B_P}(1_s)i_X(z)i_X(a^*)^* \\
& = \begin{cases} i_X(z)i_{B_P}(1_{p(z)^{-1}(p(z)\vee s)})i_X(a^*)^*
    & \text{if $p(z)\vee s < \infty$} \\
  0 & \text{otherwise} \end{cases} \\
& \in\Bb.
\end{align*}
Hence to establish \eqref{eq:span bpp}, it remains only to show that $\Bb$ is closed
under multiplication.  But Proposition~\ref{prop:yz} shows that the product
\[
i_X(x)i_{B_P}(1_s)i_X(y)^*i_X(z)i_{B_P}(1_t)i_X(w)^*
\]
of two typical generators of $\Bb$ is contained in the closed linear span of
elements of the form
\[
i_X(x)i_{B_P}(1_s)i_X(f)i_X(g)^*i_{B_P}(1_t)i_X(w)^*,
\]
which by \eqref{eq:handy2} simplifies to
\[
i_X(xf)i_{B_P}(1_{p(f)^{-1}(p(f)\vee  s) \vee p(g)^{-1}(p(g)\vee  t)})i_X(wg)^* \in\Bb.
\]

Suppose the left action of $A$ on each $X_s$ is by compact operators;
that is, $\phi_s(A) \subseteq\Kk(X_s)$ for all $s\in P$.
Let $x\in X$ and $s\in P$.  Since $X_{p(x)}$ is essential,
we can express $x = \phi_{p(x)}(a)z$ for some $a\in A$ and $z\in X_{p(x)}$.
With $\psi := \sigma \circ i_X$, we then have
\begin{align*}
\sigma(i_{B_P}(1_s)i_X(x))
& = L^\psi(1_s)\psi(x)
  = \rho^\psi_s(1)\psi_e(a)\psi(z) & \\
& = \psi^{(s)}(\phi_s(a))\psi(z) & \text{(Lemma~\ref{lemma:rho}(3))} \\
&  = \sigma(i_X^{(s)}(\phi_s(a))i_X(z)), &
\end{align*}
so $i_{B_P}(1_s)i_X(x)  = i_X^{(s)}(\phi_s(a))i_X(z) \in \Tt_{\cv}(X)$.
Since elements of the form $i_{B_P}(1_s)i_X(x)$ generate $\bpp$,
this gives $\bpp = \Tt_{\cv}(X)$.

If in addition every $s,t\in P$ have a common upper bound,
then by Proposition~\ref{prop:CP Nica} the universal map
$j_X:X\to\Oo_X$ is Nica covariant; the integrated form
$(j_X)_*:\Tt_{\cv}(X)\to\Oo_X$ is surjective
since it maps generators to generators. 
\end{proof}

\section{Faithful representations}
\label{section:faithful}

Our strategy for characterising faithful representations of $\bpp$
follows \cite[Section~5]{fowrae}.
First we use the dual coaction $\delta$ of $G$
on $\bpp$ and the canonical trace $\rho$ on $C^*(G)$
to define a positive linear map $E_\delta := (\id\otimes\rho)\circ\delta$
of norm one of $\bpp$ onto the fixed-point algebra $(\bpp)^\delta$.
When $X$ is compactly aligned,
$(B_P,P,\tau,X)$ satisfies the spanning condition
\eqref{eq:span bpp}, and $E_\delta$ is determined by
\begin{equation}\label{eq:Edelta}
E_\delta(i_X(x)i_{B_P}(1_s)i_X(y)^*) = \begin{cases}
i_X(x)i_{B_P}(1_s)i_X(y)^* & \text{if $p(x) = p(y)$} \\
0 & \text{otherwise.}
\end{cases}
\end{equation}

\begin{definition}\label{defn:amenable}
The system $(B_P,P,\tau,X)$ is {\em amenable\/} if $E_\delta$
is faithful on positive elements.
\end{definition}

The argument of \cite[Lemma~6.5]{lacarae} shows that if
$G$ is an amenable group, then the system $(B_P,P,\tau,X)$ is amenable.
In Corollary~\ref{cor:amenability} we will show that
$(B_P,P,\tau,X)$ is also amenable when $X$ is compactly aligned
and $G$ is a free product $*(G^\lambda,P^\lambda)$ with each $G^\lambda$
an amenable group.

\begin{theorem}\label{theorem:faithfulness of representations}
Suppose $(G,P)$ is a quasi-lattice ordered group
and $X$ is a compactly-aligned product system over $P$
of essential Hilbert $A$--$A$ bimodules
such that the system $(B_P,P,\tau,X)$ is amenable.
Let $\psi$ be a Nica-covariant Toeplitz representation of $X$
on a Hilbert space $\Hh$.
Then $L^\psi\times\psi$ is a faithful representation of $\bpp$
if and only if
\begin{multline}\label{eq:psifaithful}
\text{for every $n\ge 1$ and $s_1, \dots, s_n\in P\setminus\{e\}$,
  the subrepresentation} \\
a\in A\mapsto\psi_e(a)\prod_{k=1}^n \bigl(1 - L^\psi(1_{s_k})\bigr)
  \text{ of $\psi_e$ is faithful.}
\end{multline}
\end{theorem}

\begin{proof}[Proof of necessity of \eqref{eq:psifaithful}]
Let $\pi:A\to B(\Hh)$ be a faithful nondegenerate representation of $A$
on a Hilbert space $\Hh$,
let $l:X\to\Ll(F(X))$ be the Fock representation of $X$,
and let $\Psi := F(X)\dashind_A^{\Ll(X)}\pi\circ l$;
by Lemma~\ref{lemma:fock covariant}, $\Psi$ is a Nica-covariant
Toeplitz representation of $X$ on $F(X)\otimes_A\Hh$.
We claim that
\[
a\in A \mapsto \Psi_e(a)\prod_{k=1}^n \bigl(1 - L^\Psi(1_{s_k})\bigr)
\]
is faithful.
Since $L^\Psi(1_{s_k}) = \alpha^\Psi_{s_k}(1)$ is the orthogonal
projection of $F(X)\otimes_A\Hh$ onto $\bigoplus_{t\in s_kP} X_t\otimes_A\Hh$
(see the proof of Lemma~\ref{lemma:fock covariant}),
each projection $1 - L^\Psi(1_{s_k})$ dominates the projection $Q_e$ onto
the $\Psi_e$-invariant subspace $X_e \otimes_A\Hh$.
To establish the claim it thus suffices to show that
the subrepresentation $Q_e\Psi_e$ of $\Psi_e$ is faithful.
But $\Psi_e = F(X)\dashind_A^A\pi$ decomposes as
$\bigoplus_{t\in P} X_t\dashind_A^A\pi$,
so $Q_e\Psi_e = A\dashind_A^A\pi$ is unitarily equivalent to $\pi$,
and hence faithful.

Now suppose that $L^\psi\times\psi$ is faithful and $a\in A$.
Let
\[
T :=  i_{B_P}\Bigl( \prod_{k=1}^n (1  - 1_{s_k}) \Bigr)  i_X(a) \in \bpp.
\]
Then
\begin{multline*}
\norm a
= \norm{\Psi_e(a)\prod_{k=1}^n \bigl(1 - L^\Psi(1_{s_k})\bigr)}
= \norm{L^\Psi\times \Psi(T)}
\le \norm T \\
= \norm{L^\psi\times\psi(T)}
= \norm{\psi_e(a)\prod_{k=1}^n \bigl(1 - L^\psi(1_{s_k})\bigr)}
\le \norm a,
\end{multline*}
giving \eqref{eq:psifaithful}.
\end{proof}

Our proof that \eqref{eq:psifaithful} implies faithfulness
of $L^\psi\times\psi$ is based on the argument of \cite[Section~6]{fowrae}:
in Proposition~\ref{prop:fpa and EPsi}(1)
we prove that $L^\psi\times\psi$ is faithful on $(\bpp)^\delta$,
and in Proposition~\ref{prop:fpa and EPsi}(2)
we construct a spatial version $E_\psi$ of $E_\delta$ such that
$(L^\psi\times\psi)\circ E_\delta = E_\psi\circ(L^\psi\times\psi)$.
Faithfulness of $L^\psi\times\psi$ then follows easily:
if $L^\psi\times\psi(b) = 0$, then
\[
0 = E_\psi\circ(L^\psi\times\psi)(b^*b) = (L^\psi\times\psi)\circ E_\delta(b^*b),
\]
so by Proposition~\ref{prop:fpa and EPsi}(1), $E_\delta(b^*b) = 0$.
The amenability hypothesis then forces $b^*b = 0$, and hence $b=0$. 

We begin by reviewing some notation and results from \cite[Remark~1.5]{lacarae}
and \cite[Remark~5.2]{fowrae}.
Let $F$ be a finite subset of $P$.  A subset $C$ of $F$ is an
{\em initial segment\/} of $F$ if $c := \sigma C$ is finite
and $C = \{t\in F: t\le c\}$.
(Recall that $\sigma C$ is the least upper bound of $C$;
we use the convention that $\sigma\emptyset = e$.)
For each such $C$ there is a nonzero projection $Q_C$ in $B_P$ defined by
\[
Q_C := 1_c \prod_{\{t\in F: c < t\vee c < \infty\}} (1 - 1_t),
\]
and as $C$ ranges over the initial segements of $F$, these projections
form a decomposition of the identity in $B_P$.

\begin{lemma}\label{lemma:QC}
Suppose $(G,P)$ is a quasi-lattice ordered group,
$X$ is a product system over $P$ of essential Hilbert $A$--$A$ bimodules,
$\psi$ is a Nica-covariant Toeplitz representation of $X$ on $\Hh$,
$F$ is a finite subset of $P$,
$C$ is an initial segment of $F$,
$x,y\in X$ and $s\in P$.
Let $c = \sigma C$, so that $C = \{t \in F: t \le c\}$.

\textup{(1)} If $p(x) = p(y)$,
then the operator $\psi(x)L^\psi(1_s)\psi(y)^*$
is in the commutant of $L^\psi(B_P)$.
In particular, it commutes with $L^\psi(Q_C)$.

\textup{(2)} If $p(x)s$, $p(y)s\in F$, then
\begin{multline*}
L^\psi(Q_C)\psi(x)L^\psi(1_s)\psi(y)^*L^\psi(Q_C) \\
  = \left\{ \begin{array}{l}
      L^\psi(Q_C)
        \psi(x)L^\psi(1_{p(x)^{-1}c})L^\psi(1_{p(y)^{-1}c})\psi(y)^*
        L^\psi(Q_C) \\
      \phantom{0} \qquad\qquad\qquad\qquad \text{if $p(x)s\le c$ and $p(y)s\le c$} \\
      0 \qquad \text{otherwise.}
\end{array}
\right.
\end{multline*}
\end{lemma}

\begin{proof} The proof, based on Lemma~\ref{lemma:handy},
is identical in form to the proof of \cite[Lemma~5.3]{fowrae}.
\end{proof}

\begin{lemma}\label{lemma:compute norm}
Suppose $(G,P)$ is a quasi-lattice ordered group,
$X$ is a product system over $P$ of essential Hilbert $A$--$A$ bimodules,
and $\psi$ is a Nica-covariant Toeplitz representation of $X$ which satisfies
\eqref{eq:psifaithful}.
Suppose further that $F$ be a finite subset of $P$
and $Z$ is a finite sum $\sum \psi(x_k)L^\psi(1_{s_k})\psi(y_k)^*$
such that $p(x_k)s_k = p(y_k)s_k \in F$ for each $k$.
Then
\begin{equation}\label{eq:normTC}
\norm Z = \max\{\norm{T_C}: \text{$C$ is an initial segment of $F$}\},
\end{equation}
where $T_C$ is the adjointable operator on $X_{\sigma C}$ defined by
\begin{equation}\label{eq:TC}
T_C := \sum_{p(x_k)s_k \le \sigma C}
       \Theta_{x_k,y_k} \otimes_A 1^{p(x_k)^{-1}\sigma C}
\end{equation}
\end{lemma}

\begin{proof}
Since $\{Q_C: \text{$C$ is an initial segment of $F$}\}$
is a decomposition of the identity in $B_P$,
and since $L^\psi$ is a unital representation of $B_P$,
the projections $L^\psi(Q_C)$ decompose the identity operator.
By Lemma~\ref{lemma:QC}(1), $Z$ commutes with each $L^\psi(Q_C)$, and thus
\[
\norm{Z} = \max\{\norm{L^\psi(Q_C)Z}:
\text{$C$ is an initial segment of $F$.}\}
\]
Fix an initial segment $C$, and let $c:=\sigma C$.
By Lemma~\ref{lemma:QC}(2) and Lemma~\ref{lemma:rho}(6),
\begin{align*}
L^\psi(Q_C)Z
& = L^\psi(Q_C) \sum \psi(x_k)L^\psi(1_{s_k})\psi(y_k)^* \\
& = L^\psi(Q_C) \sum_{p(x_k)s_k\le c}
                      \psi(x_k)L^\psi(1_{p(x_k)^{-1}c})\psi(y_k)^* \\
& = L^\psi(Q_C) \sum_{p(x_k)s_k\le c}
                      \rho^\psi_c ( \Theta_{x_k,y_k} \otimes_A 1) \\
& = L^\psi(Q_C) \rho^\psi_c(T_C),
\end{align*}
so it suffices to show that
\begin{equation}\label{eq:compute norms2}
\norm{L^\psi(Q_C) \rho^\psi_c(T_C)} = \norm{T_C}.
\end{equation}

Let
\begin{equation}\label{eq:RC}
R_C := \prod_{\{t\in F: c < t \vee c < \infty\}}
     ( 1 - 1_{c^{-1}(t \vee c)}) \in B_P.
\end{equation}
Since $\psi$ satisfies \eqref{eq:psifaithful},
\[
a \mapsto \psi_e(a)\prod_{\{t\in F: c < t \vee c < \infty\}}
     ( 1 - L^\psi(1_{c^{-1}(t \vee c)}) )
= \psi_e(a)L^\psi(R_C)
\]
is a faithful representation of $A$.
By Lemma~\ref{lemma:rho}(4), the representation
$T\in\Ll(X_c) \mapsto \alpha^\psi_c(L^\psi(R_C))\rho^\psi_c(T)$
is thus also faithful.
But $\alpha^\psi_c(L^\psi(R_C)) = L^\psi(\tau_c(R_C)) = L^\psi(Q_C)$,
and hence \eqref{eq:compute norms2} is satisfied.
\end{proof}

\begin{prop}\label{prop:fpa and EPsi}
Suppose $(G,P)$ is a quasi-lattice ordered group,
$X$ is a compactly-aligned product system over $P$
of essential Hilbert $A$--$A$ bimodules,
and $\psi$ is a Nica-covariant Toeplitz representation of $X$ which satisfies
\eqref{eq:psifaithful}.

\textup{(1)} $L^\psi\times\psi$ is isometric on $(\bpp)^\delta$.

\textup{(2)} There is a linear map $E_\psi$ of norm one of
$L^\psi\times\psi(\bpp)$ onto $L^\psi\times\psi\bigl((\bpp)^\delta\bigr)$
such that $E_\psi \circ (L^\psi\times\psi) = (L^\psi\times\psi) \circ E_\delta$.
\end{prop}

\begin{proof}
(1) Since $X$ is compactly aligned, the spanning condition
\eqref{eq:span bpp} holds.  Since $E_\delta$ is continuous
and maps onto $(\bpp)^\delta$,
we deduce that finite sums
\[
z := \sum i_X(x_k)i_{B_P}(1_{s_k})i_X(y_k)^*
\]
in which $p(x_k) = p(y_k)$ for all $k$ are dense in $(\bpp)^\delta$.
It therefore suffices to fix such a $z$
and show that $\norm{L^\psi\times\psi(z)} = \norm z$.

Let $\sigma$ be a faithful nondegenerate representation
of $\bpp$ such that $(\overline\sigma\circ i_{B_P},\sigma\circ i_X)$
is a covariant representation of $(B_P,P,\tau,X)$.
By Proposition~\ref{prop:Lpsi}, $i := \sigma\circ i_X$ is a covariant
representation of $X$ and $\overline\sigma\circ i_{B_P} = L^i$.
Since $L^i\times i = \sigma$ is faithful, $i$ satisfies \eqref{eq:psifaithful}.
Hence with $F := \{p(x_k)s_k\}$,
Lemma~\ref{lemma:compute norm} gives
\begin{align*}
\norm{L^\psi\times\psi(z)}
& = \norm{\sum \psi(x_k)L^\psi(1_{s_k})\psi(y_k)^*} \\
& = \max\{\norm{T_C}: \text{ $C$ is an initial segment of $F$}\} \\
& = \norm{\sum i(x_k)L^i(1_{s_k})i(y_k)^*}
= \norm{L^i\times i(z)}
= \norm z.
\end{align*}

(2) Since $X$ is compactly aligned,
finite sums of the form
\[
w := \sum i_X(x_k)i_{B_P}(1_{s_k})i_X(y_k)^*
\]
are dense in $\bpp$.
We will show that
$\norm{L^\psi\times\psi(E_\delta(w))} \le \norm{L^\psi\times\psi(w)}$;
it follows that $E_\psi$
is well-defined on operators of the form $L^\psi\times\psi(w)$
and extends to the desired linear contraction.

Let $F := \{p(x_k)s_k\} \cup \{p(y_k)s_k\}$,
and let $Z := L^\psi\times\psi(E_\delta(w))$;
by \eqref{eq:Edelta},
\[
Z = \sum_{p(x_k) = p(y_k)} \psi(x_k)L^\psi(1_{s_k})\psi(y_k)^*.
\]
By Lemma~\ref{lemma:compute norm}, there is an initial
segment $C$ of $F$ such that
$\norm Z = \norm{T_C}$.
Let $c := \sigma C$.  We will construct a projection $R\in B_P$ such that
$a\in A \mapsto \psi_e(a)L^\psi(R)$ is faithful,
then define $Q := L^\psi(\tau_c(R)) = \alpha^\psi_c(L^\psi(R))$,
and show that $Q(L^\psi\times\psi(w))Q = Q\rho^\psi_c(T_C)$.
This will complete
the proof, since by Lemma~\ref{lemma:rho}(4) we then have
\[
\norm Z = \norm{T_C} = \norm{Q\rho^\psi_c(T_C)} = 
\norm{Q(L^\psi\times\psi(w))Q} \le \norm{L^\psi\times\psi(w)}.
\]

For each $s,t\in C$ such that $s \ne t$ and
$s^{-1}c \vee t^{-1}c < \infty$,
define $d_{s,t}\in P$ as in \cite[Lemma~3.2]{lacarae}:
\[
d_{s,t} =
\begin{cases}
  (s^{-1}c)^{-1}(s^{-1}c \vee t^{-1}c)
    & \text{if $s^{-1}c < s^{-1}c \vee t^{-1}c$} \\
  (t^{-1}c)^{-1}(s^{-1}c \vee t^{-1}c)
    & \text{otherwise,} \\
\end{cases}
\]
noting in particular that $d_{s,t}$ is never the identity in $P$.
Let $R_C$ be as in \eqref{eq:RC}, and define
\[
R := R_C \prod_{\substack{s\ne t\in C\\ s^{-1}c \vee t^{-1}c < \infty}}
( 1 - 1_{d_{s,t}}).
\]
By condition \eqref{eq:psifaithful},
$a\in A \mapsto L^\psi(R)\psi_e(a)$ is faithful.
The proof that $Q(L^\psi\times\psi(w))Q = Q\rho^\psi_c(T_C)$ is exactly as in
\cite[Proposition~5.5]{fowrae}, so we omit it.
\end{proof}

\begin{prop}\label{prop:EPsi faithful}
Suppose $(G,P)$ is a quasi-lattice ordered group
and $X$ is a compactly-aligned product system over $P$
of essential Hilbert $A$--$A$ bimodules.
Let $\pi$ be a nondegenerate representation of $A$ on a Hilbert space $\Hh$,
and let $\Psi$ be the representation $F(X)\dashind_A^{\Ll(X)}\pi\circ l$,
where $l:X\to\Ll(F(X))$ is the Fock representation of $X$.
There is a projection $E_\Psi$ of norm one of $L^\Psi\times \Psi(\bpp)$ onto
$L^\Psi\times \Psi((\bpp)^\delta)$ such that
\begin{equation}\label{eq:EPsi}
E_\Psi \circ (L^\Psi\times\Psi) = (L^\Psi\times\Psi) \circ E_\delta;
\end{equation}
moreover,  $E_\Psi$ is faithful on positive operators.
\end{prop}

\begin{proof} Denote by $Q_t$ the orthogonal projection
of $F(X)\otimes_A\Hh$ onto $X_t\otimes_A\Hh$.
Since the $Q_t$'s are mutually orthogonal, the formula
\[
E_\Psi(T) := \sum_{t\in P} Q_t T Q_t
\qquad\text{for $T\in L^\Psi\times \Psi(\bpp)$}
\]
defines a completely positive projection of norm one
which is faithful on positive operators.
We claim that
\begin{equation}\label{eq:EPsi2}
E_\Psi(\Psi(x)L^\Psi(1_s)\Psi(y)^*)
= \begin{cases}
\Psi(x)L^\Psi(1_s)\Psi(y)^* & \text{if $p(x) = p(y)$} \\
0 & \text{otherwise.}
\end{cases}
\end{equation}
Since $X$ is compactly aligned the spanning condition
\eqref{eq:span bpp} holds, and hence \eqref{eq:EPsi}
follows from \eqref{eq:EPsi2} and \eqref{eq:Edelta}.

Suppose $x,y\in X$ and $s\in P$.
For each $t \in P$, $\Psi(x)L^\Psi(1_s)\Psi(y)^*$ is zero on $X_t\otimes_A\Hh$
unless $p(y)s \le t$, in which case $\Psi(x)L^\Psi(1_s)\Psi(y)^*$ maps
$X_t\otimes_A\Hh$ into $X_{p(x)p(y)^{-1}t}\otimes_A\Hh$.
Thus if $p(x) \ne p(y)$, $Q_t\Psi(x)L^\Psi(1_s)\Psi(y)^*Q_t=0$ for every $t\in P$,
and $E_\Psi(\Psi(x)L^\Psi(1_s)\Psi(y)^*) = 0$.
If on the other hand $p(x) = p(y)$, then
$Q_t\Psi(x)L^\Psi(1_s)\Psi(y)^*Q_t = \Psi(x)L^\Psi(1_s)\Psi(y)^*Q_t$
for each $t\in P$, and thus
\begin{multline*}
E_\Psi(\Psi(x)L^\Psi(1_s)\Psi(y)^*)
= \sum_{t\in P} Q_t\Psi(x)L^\Psi(1_s)\Psi(y)^*Q_t \\
= \Psi(x)L^\Psi(1_s)\Psi(y)^*\sum_{t\in P} Q_t
= \Psi(x)L^\Psi(1_s)\Psi(y)^*.
\end{multline*}
\end{proof}

\begin{cor}\label{cor:EPsi faithful}
Suppose $\pi$ is faithful.  Then the system $(B_P,P,\tau,X)$ is amenable
if and only if the representation $L^\Psi\times\Psi$ of $\bpp$ is faithful.
\end{cor}

\begin{proof} Suppose $L^\Psi\times\Psi$ is faithful.
By Proposition~\ref{prop:EPsi faithful},
$(L^\Psi\times\Psi) \circ E_\delta = E_\Psi \circ (L^\Psi\times\Psi)$ is faithful
on positive elements, hence so is $E_\delta$;
that is, $(B_P,P,\tau,X)$ is amenable.
Since $\Psi$ satisfies \eqref{eq:psifaithful}
(see the proof of necessity of \eqref{eq:psifaithful}),
the converse follows from Theorem~\ref{theorem:faithfulness of representations}.
\end{proof}

\section{Amenability}\label{section:amenability}

\begin{theorem}\label{theorem:amenability}
Suppose $\theta: (G, P) \to (\Gg, \Pp)$ is a homomorphism of quasi-lattice
ordered groups such that, whenever $s\vee t<\infty$,
\begin{equation}\label{eq:theta}
  \theta(s \vee t)=  \theta(s) \vee \theta(t)\ \text{ and }\ 
  \theta(s)=  \theta(t)\Longrightarrow s = t,
\end{equation}
and suppose that $\Gg$ is amenable.
If $X$ is a compactly-aligned product system over $P$
of essential Hilbert $A$--$A$ bimodules,
then the system $(B_P,P,\tau,X)$ is amenable.
\end{theorem}

\begin{proof}
Our proof is essentially that of \cite[Theorem~6.1]{fowrae},
suitably modified to handle Hilbert bimodules.
The homomorphism $\theta:G \to \Gg$ induces a coaction
$\delta_\theta = (\id \otimes \theta) \circ \delta$ of $\Gg$ on $\bpp$,
and hence a conditional expectation $E_{\delta_\theta}$
of $\bpp$ onto the fixed-point algebra $(\bpp)^{\delta_\theta}$,
such that 
\[
E_{\delta_\theta}(i_X(x)i_{B_P}(1_s)i_X(y)^*) =
\begin{cases}
  i_X(x)i_{B_P}(1_s)i_X(y)^*& \text{if $\theta(p(x)) = \theta(p(y))$} \\
  0 & \text{otherwise.}
\end{cases}
\]
Since $\Gg$ is amenable, $E_{\delta_\theta}$ is faithful on positive elements.

Let $l:X\to\Ll(F(X))$ be the Fock representation of $X$,
let $\pi$ be a faithful nondegenerate representation of $A$
on a Hilbert space $\Hh$, and let $\Psi := F(X)\dashind_A^{\Ll(X)}\pi\circ l$.
By Proposition~\ref{prop:EPsi faithful}, for every $b\in\bpp$ we have
\[
(L^\Psi\times\Psi)\circ E_\delta(b)
= E_\Psi(L^\Psi\times\Psi(E_{\delta_\theta}(b))).
\]
Since $E_{\delta_\theta}$ and $E_\Psi$ are faithful on positive elements,
to show that $(B_P,P,\tau,X)$ is amenable it suffices to show that
$L^\Psi\times\Psi$ is faithful on $(\bpp)^{\delta_\theta}$.

Let $\sigma$ be a faithful representation of
$\bpp$ such that $(\overline\sigma \circ i_{B_P}, \sigma \circ i_X)$
is a covariant representation of $(B_P, P, \tau, X)$.
By Proposition~\ref{prop:Lpsi},
$i = \sigma \circ i_X$ is a covariant representation of $X$
and $\overline\sigma \circ i_{B_P} = L^i$.
Observe that $i$ is isometric since, by Lemma~\ref{lemma:fock covariant},
\begin{align*}
\norm x
= \norm{\Psi(x)}
& = \norm{(L^\Psi\times \Psi)\circ i_X(x)} \\
& \le \norm{i_X(x)}
= \norm{\sigma\circ i_X(x)}
= \norm{i(x)}
\le \norm x.
\end{align*}

Let $\Ff$ be the set of all finite subsets $F$ of $\Pp$
which are closed under $\vee$ in the sense that $s\vee t\in F$
whenever $s,t\in F$ and $s\vee t < \infty$.
Exactly as in the proof of \cite[Theorem~6.1]{fowrae},
one can use Proposition~\ref{prop:yz} to show that,
for each $F\in\Ff$,
\[
\Uu_F := \clsp\{i_X(x)i_{B_P}(1_s)i_X(y)^*: \theta(p(x)s) = \theta(p(y)s) \in F\}
\]
is a $C^*$\ndash subalgebra of $\bpp$.
Applying $\Phi_{\delta_\theta}$ to both sides of \eqref{eq:span bpp} gives
\[
(\bpp)^{\delta_\theta}
  = \clsp\{i_X(x)i_{B_P}(1_s)i_X(y)^*: \theta(p(x)) = \theta(p(y))\};
\]
since $\Ff$ is directed under set inclusion
(see the proof of \cite[Lemma~4.1]{lacarae}),
we deduce that
\[
(\bpp)^{\delta_\theta} = \overline{\bigcup_{F\in\Ff} \Uu_F}.
\]
By \cite[Lemma~1.3]{alnr}, to prove that $L^\Psi\times\Psi$ is faithful
on $(\bpp)^{\delta_\theta}$ it is enough to prove it is faithful
on each of the subalgebras $\Uu_F$.
We shall accomplish this by inducting on $\lvert F \rvert$.

First suppose $F = \{r\}$ for some $r\in\Pp$.
Let $W_r$ be the Hilbert $A$--$A$ bimodule
$\bigoplus_{t\in\theta^{-1}(r)} X_t$.
We claim that, for each Nica-covariant Toeplitz representation $\psi$ of
$X$ on a Hilbert space $\Kk$,
there is a linear map $\psi_r:W_r\to B(\Kk)$ 
which satisfies
$\psi_r(\oplus x_t) = \sum \psi_t(x_t)$,
and that $(\psi_r,\psi_e)$ is then a Toeplitz representation of $W_r$.
First observe that if $x,y\in X$ satisfy
$p(x) \ne p(y)$ and $\theta(p(x)) = \theta(p(y)) = r$,
then by \eqref{eq:theta} we have $p(x) \vee p(y) = \infty$,
and hence $\psi(x)^*\psi(y) = 0$.
Now suppose $\oplus x_t$ belongs to the algebraic direct sum
$\bigodot_{t\in\theta^{-1}(r)} X_t$; such vectors are dense in $W_r$.
Then
\begin{align*}
\norm{\sum_t \psi_t(x_t)}^2
& = \norm{\sum_{t,t'} \psi_t(x_t)^*\psi_{t'}(x_{t'})}
= \norm{\sum_t \psi_t(x_t)^*\psi_t(x_t)} \\
& = \norm{\sum_t \psi_e(\langle x_t,x_t \rangle_A)}
\le \norm{\sum_t \langle x_t, x_t \rangle_A} \\
& = \norm{\langle \oplus x_t, \oplus x_t \rangle_A}
= \norm{\oplus x_t}^2,
\end{align*}
ensuring the existence of $\psi_r$.
It is  routine to check that 
$(\psi_r,\psi_e)$ is a Toeplitz representation of $W_r$.
Write $\alpha^\psi_r$ for the endomorphism of $\psi_e(A)'$
which corresponds to $(\psi_r,\psi_e)$
(Proposition~\ref{prop:alpha}),
and write $\rho^\psi_r$ for the associated representation of $\Ll(W_r)$
(Lemma~\ref{lemma:rho}).

Suppose $Z$ is a finite sum
$\sum i_X(x_k)i_{B_P}(1_{s_k})i_X(y_k)^*$
such that
$\theta(p(x_k)s_k) = \theta(p(y_k)s_k) = r$
for every $k$;
to prove $L^\Psi\times\Psi$ faithful on $\Uu_{\{r\}}$
we will show that $\norm{L^\Psi\times\Psi(Z)} = \norm Z$.
For each $k$, let $\Theta_{x_k,y_k}\otimes_A 1^{s_k}$
denote the operator in $\Ll(W_r)$  which is the image of
\[
\Theta_{x_k,y_k}\in\Kk(X_{p(y_k)},X_{p(x_k)})
\mapsto \Theta_{x_k,y_k}\otimes_A 1^{s_k}
\in \Ll(X_{p(y_k)s_k},X_{p(x_k)s_k})
\subset \Ll(W_r).
\]
Define $T := \sum \Theta_{x_k,y_k}\otimes_A 1^{s_k} \in \Ll(W_r)$.
It is routine to check that
\[
\rho^\Psi_r(T)
= \sum \Psi(x_k)L^\Psi(1_{s_k})\Psi(y_k)^*
= L^\Psi\times\Psi(Z),
\]
and similarly $\rho^i_r(T) = L^i \times i(Z) = \sigma(Z)$.
Since $\Psi_e$ and $i_e$ are faithful representations of $A$,
the representations $\rho^\Psi_r$ and $\rho^i_r$ are isometric, and thus
\[
\norm{L^\Psi\times\Psi(Z)} = \norm{\rho^\Psi_r(T)} = \norm T
  = \norm{\rho^i_r(T)} = \norm{\sigma(Z)} = \norm Z.
\]

For the inductive step,
suppose $F \in \Ff$ and $L^\Psi\times\Psi$ is faithful on $\Uu_{F'}$
whenever $F' \in \Ff$ and $\lvert F' \rvert < \lvert F \rvert$;
we aim to prove that $L^\Psi\times\Psi$ is faithful on $\Uu_F$. 
Since $F$ is finite it has a minimal element;
that is, there exists $r_0 \in F$ such that
$r_0 < r_0 \vee r$ for each $r\in F \setminus\{r_0\}$.
As in the proof of \cite[Theorem~6.1]{fowrae}
we have $L^\Psi\times\Psi(\Uu_{\{r\}})P_{r_0} = \{0\}$ for each
$r \in F \setminus\{r_0\}$,
where $P_{r_0}$ denotes the orthogonal projection of
$F(X) \otimes_A \Hh$
onto $\bigoplus_{t \in \theta^{-1}(r_0)} X_t\otimes_A\Hh$.

On the other hand, we have already demonstrated that
$L^\Psi\times\Psi$ maps $\Uu_{r_0}$ isometrically into
the range of $\rho^\Psi_{r_0}$, and an easy calculation
shows that $P_{r_0} = \alpha^\Psi_{r_0}(Q_e)$,
where $Q_e$ is the orthogonal projection onto $X_e\otimes_A\Hh$.
Since $a\mapsto\Psi_e(a)Q_e$ is faithful,
by Lemma~\ref{lemma:rho}(4) the representation
$S\in\Ll(W_{r_0}) \mapsto P_{r_0}\rho^\Psi_{r_0}(S)$ is also faithful.
Hence the map $Y \in \Uu_{r_0} \mapsto L^\Psi\times\Psi(Y)P_{r_0}$
is faithful.

Now suppose $Y \in \Uu_F$ and $L^\Psi\times\Psi(Y) = 0$.
We will show that $Y \in \Uu_{F \setminus \{r_0\}}$,
from which the inductive hypothesis implies that $Y = 0$.
Let $(Y_n)$ be a sequence in
\[
\Span\{i_X(x)i_{B_P}(1_s)i_X(y)^*: \theta(p(x)s) = \theta(p(y)s) \in F\}
\]
which converges in norm to $Y$, and express each $Y_n$ as a sum
$\sum_{r\in F} Y_{n,r}$, where $Y_{n,r} \in \Uu_{\{r\}}$.
For each $n$,
\[
\norm{L^\Psi\times\Psi(Y_n) P_{r_0}}
 = \norm{L^\Psi\times\Psi(Y_{n,r_0}) P_{r_0}}
 = \norm{Y_{n,r_0}},
\]
and consequently $Y_{n, r_0} \to 0$.
Thus $Y_n - Y_{n,r_0} \to Y$, which shows that
$Y \in \Uu_{F \setminus\{r_0\}}$, as claimed.
\end{proof}

\begin{cor}\label{cor:amenability}
Suppose $(G^\lambda, P^\lambda)$ is a quasi-lattice ordered
group with $G^\lambda$ am\-en\-able for each $\lambda$ belonging
to some index set $\Lambda$.
If $X$ is a compactly-aligned product system over $P := * P^\lambda$,
then the system $(B_P,P,\tau,X)$ is amenable.
\end{cor}

\begin{proof} The group $\bigoplus G^\lambda$ is amenable,
and by \cite[Proposition~4.3]{lacarae} the canonical map
$\theta: * G^\lambda \to \bigoplus G^\lambda$ satisfies \eqref{eq:theta}. 
\end{proof}

\section{Applications}\label{section:applications}

In Section~\ref{section:crossed products},
we associated with each twisted semigroup dynamical system
$(A,P,\beta,\omega)$ a product system $X = X(A,P,\beta,\omega)$
of essential Hilbert $A$--$A$ bimodules over the opposite semigroup $P^o$
(Lemma~\ref{lemma:crossed product ps}),
and we showed that the Cuntz-Pimsner algebra $\Oo_X$ is canonically isomorphic
to the crossed product $A\cross_{\beta,\omega} P$;
we also showed that  $\Tt_X$ has the structure of a certain
``Toeplitz'' crossed product $\Tt(A\cross_{\beta,\omega} P)$
(Proposition~\ref{prop:crossed product}).
Suppose now that $(G^o,P^o)$ is quasi-lattice ordered;
this is equivalent to $(G,P)$ being quasi-latticed ordered in its
{\em right\/}-invariant partial order
($s\le t \Leftrightarrow ts^{-1} \in P$).
Since the left action of $A$ on each fiber $X_s$ is by compact operators,
$X$ is compactly aligned (Lemma~\ref{prop:CA conditions})
and $\Tt_{\cv}(X) = \bpp$ (Theorem~\ref{theorem:subalgebra}).
Hence we can apply Theorem~\ref{theorem:faithfulness of representations}
to characterize the faithful representations of $\Tt_{\cv}(X)$.
This is particularly helpful when $(G^o,P^o)$ is a total order
since $\Tt_{\cv}(X) = \Tt_X$; more generally,
when every $s,t\in P^o$ have a common upper bound in $P^o$
(i.e. $Ps\cap Pt \ne \emptyset$),
the crossed product $A\cross_{\beta,\omega} P = \Oo_X$ is a quotient
of $\Tt_{\cv}(X)$ (Theorem~\ref{theorem:subalgebra}).

We begin by showing that $\Tt_{\cv}(X)$, too, has a crossed product structure:

\begin{definition}
Suppose $P$ is a subsemigroup of a group $G$
and $(G^o,P^o)$ is quasi-lattice ordered.
A {\em Nica-Toeplitz covariant representation\/} of $(A,P,\beta,\omega)$
is a Toeplitz covariant representation $(\pi,V)$
such that
\begin{equation}\label{eq:V covariant}
V_s^*V_sV_t^*V_t =
\begin{cases} V_{s\vee t}^*V_{s\vee t}
  & \text{if $s\vee t < \infty$} \\
0 & \text{otherwise,} \end{cases}
\end{equation}
where $s\vee t$ denotes the least upper bound of $s$ and $t$
in the right-invariant partial order on $(G,P)$.
\end{definition}

The following Proposition establishes the existence
of a $C^*$\ndash al\-ge\-bra which is universal for such pairs $(\pi,V)$,
as in Definition~\ref{defn:omega crossed product}.
We call this algebra the {\em Nica-Toeplitz crossed product\/}
of $(A,P,\beta,\omega)$, and denote it
$\Tt_{\cv}(A\cross_{\beta,\omega} P)$.
Let $i_X:X\to\Tt_{\cv}(X)$ be universal for Nica-covariant
Toeplitz representations of $X$.
Lemma~\ref{lemma:strict convergence} is easily adapted to this setting,
and allows us to define $i_P:P\to M\Tt_{\cv}(X)$ by
$i_P(s) = \lim i_X(s,\beta_s(a_i))^*$;
here  $(a_i)$ is an approximate identity for $A$,
and the convergence is strict.
We also define $i_A:A\to\Tt_{\cv}(X)$ by $i_A(a) := i_X(e,a)$.

\begin{prop}\label{prop:crossed product2}
$(\Tt_{\cv}(X),i_A, i_P)$
is a Nica-Toeplitz crossed product for $(A,P,\beta,\omega)$.
\end{prop}

\begin{proof}
As in the proof of Proposition~\ref{prop:crossed product},
$i_A$ is nondegenerate.
We verify the obvious analogues of conditions (a), (b), and (c)
in Definition~\ref{defn:omega crossed product}.
For (a), let $\sigma$  be a nondegenerate representation of $\Tt_{\cv}(X)$
on a Hilbert space $\Hh$,
let $\pi := \sigma\circ i_A$, and let $V := \overline\sigma\circ i_P$;
we must show that $(\pi,V)$ is a Nica-Toeplitz covariant representation
of $(A,P,\beta,\omega)$.
Exactly as in the proof of Proposition~\ref{prop:crossed product},
$(\pi,V)$ is a Toeplitz covariant representation of $(A,P,\beta,\omega)$,
so we need to establish \eqref{eq:V covariant}.
Fix $s\in P$.
For any $a\in A$ and $h\in\Hh$ we have
\begin{align*}
V_s^*\pi(a)h
& = \sigma(i_P(s)^*i_A(a))h
= \sigma(\lim i_X(s,\beta_s(a_i))i_X(e,a))h \\
& = \sigma(\lim i_X(s,\beta_s(a_i)a))h
= \sigma \circ i_X(s, \overline{\beta_s}(1)a)h,
\end{align*}
and since $\pi$ is nondegenerate this shows that
\[
V_s^*V_s\Hh
= \clsp\{\sigma\circ i_X(\xi)h: \xi\in X_s, h\in\Hh\}
= \alpha^{\sigma\circ i_X}_s(1).
\]
Since $X$ is compactly aligned,
$\sigma\circ i_X$ is Nica covariant
(Theorem~\ref{theorem:subalgebra} and Proposition~\ref{prop:preserve covariance}),
and \eqref{eq:V covariant} follows.

For condition (b), let $(\pi,V)$ be any Nica-Toeplitz covariant representation
on $\Hh$.
As in the proof of Proposition~\ref{prop:crossed product},
$\psi(s,x) := V_s^*\pi(x)$ defines a nondegenerate Toeplitz covariant representation
$\psi:X\to B(\Hh)$.   To see that it is Nica-covariant,
let $s\in P$, and note that for  any $a\in A$ we have
\begin{align*}
\psi(s,\overline{\beta_s}(1)a)
& = \lim \psi(s,\beta_s(a_i)a)
 = \lim V_s^*\pi(\beta_s(a_i)a) \\
& = \lim V_s^*V_s\pi(a_i)V_s^*\pi(a)
 = V_s^*\pi(a).
\end{align*}
Since $\pi$ is nondegenerate, this implies  that
$\alpha^\psi_s(1) = V_s^*V_s$,
and hence $\psi$ is Nica covariant by \eqref{eq:V covariant}.
Defining $\pi\times V := \psi_*:\Tt_{\cv}(X)\to B(\Hh)$
gives the desired representation satisfying $(\pi\times V)\circ i_A = \pi$
and $\overline{\pi\times V}\circ i_P = V$.
Condition (c) is satisfied since $i_A(a)i_P(s) = i_X(s,\overline{\beta_s}(1)a^*)^*$,
and elements of this form generate $\Tt_{\cv}(X)$.
\end{proof}

Let $(G_i,P_i)$ be a collection of abelian lattice-ordered groups.
Since $(G_i,P_i)$ is quasi-lattice ordered
in both its left and its right-invariant partial order,
so is the free product $*(G_i, P_i)$.

\begin{theorem}\label{theorem:crossed products}
Suppose $(G,P) = *(G_i, P_i)$ is a free product of abelian lattice-ordered groups
and $(\pi,V)$ is a Nica-Toeplitz covariant representation
of the twisted semigroup dynamical system $(A,P,\beta,\omega)$
on a Hilbert space $\Hh$.
Then the integrated form $\pi\times V$ is a faithful representation
of  $\Tt_{\cv}(A\cross_{\beta,\omega} P)$
if and only if
\begin{multline*}
\text{for every $n\ge 1$ and $s_1, \dots, s_n\in P\setminus\{e\}$,} \\
\text{$\pi$ acts faithfully on the range of
      $\prod_{k=1}^n \bigl(1 - V_{s_k}^*V_{s_k})$.}
\end{multline*}
\end{theorem}

\begin{proof} Let $\theta$ be the canonical homomorphism of
$*(G_i,P_i)$ onto $\bigoplus (G_i,P_i)$.
By \cite[Proposition~4.3]{lacarae}, $\theta$ satisfies
the hypotheses of Theorem~\ref{theorem:amenability};
since $X = X(A,P,\beta,\omega)$ is compactly aligned,
the system $(B_P,P,\tau,X)$ is therefore amenable.
Identifying $\Tt_{\cv}(A\cross_{\beta,\omega} P)$
with $\Tt_{\cv}(X)$ as in the previous Proposition and defining
$\psi(s,x) := V_s^*\pi(x)$, the initial projection $V_s^*V_s$
is precisely $\alpha^\psi_s(1)$,
and the result follows
from Theorem~\ref{theorem:faithfulness of representations}.
\end{proof}

Nica covariance is automatic when $(G,P)$ is totally ordered:

\begin{cor}
Suppose $(G,P)$ is a totally ordered abelian group and $(\pi,V)$ is a Toeplitz
covariant representation of $(A,P,\beta,\omega)$ on a Hilbert space $\Hh$.
Then the integrated form $\pi\times V$ is a faithful representation
of  $\Tt(A\cross_{\beta,\omega} P)$
if and only if
$\pi$ acts faithfully on $(V_s^*\Hh)^\perp$ for every $s\in P\setminus\{e\}$.
\end{cor}

\begin{cor}
Suppose $\beta$ is an extendible endomorphism of $A$.
If $(\pi,V)$ is a Toeplitz representation of $(A,\NN,\beta)$,
then $\pi\times V$ is a faithful representation of
$\Tt(A\cross_\beta\NN)$ if and only if $\pi$ acts faithfully on $(V^*\Hh)^\perp$.
\end{cor}

\subsection*{Bicovariance}
Suppose $(G,P)$ is a quasi-lattice ordered group.
Following \cite{lacarae}, in \cite{fowrae} it was shown that
$B_P\cross_{\tau,\omega} P$ is universal for
isometric $\omega$\ndash representations of $P$
which are Nica covariant;
that is, which satisfy
\begin{equation}\label{eq:V covariant2}
V_sV_s^*V_tV_t^* =
\begin{cases} V_{s\vee t}V_{s\vee t}^*
  & \text{if $s\vee t < \infty$} \\
0 & \text{otherwise.} \end{cases}
\end{equation}
Assuming that $(G^o,P^o)$ is also quasi-lattice  ordered,
we  now show that the Nica-Toeplitz crossed product
$\Tt_{\cv}(B_P\cross_{\tau,\omega}P)$ is universal for
partial isometric $\omega$\ndash representations of $P$
which are {\em bicovariant\/} in that they satisfy both
\eqref{eq:V covariant2} and \eqref{eq:V covariant}.
Note that bicovariance is automatic when $(G,P)$ is a totally ordered abelian group.

\begin{prop}\label{prop:universal bicovariant}
$i_P:P\to\Tt_{\cv}(B_P\cross_{\tau,\omega}P)$
is a bicovariant partial isometric $\omega$-represent\-ation of $P$
whose range generates $\Tt_{\cv}(B_P\cross_{\tau,\omega}P)$ as a $C^*$\ndash al\-ge\-bra.
Moreover, for every bicovariant partial isometric
$\omega$-re\-pre\-senta\-tion $V$,
there is  a representation $V_*$  of $\Tt_{\cv}(B_P\cross_{\tau,\omega}P)$
such that $V_*\circ i_P = V$.
\end{prop}

\begin{proof}
Let $\sigma$ be a faithful nondegenerate representation of
$\Tt_{\cv}(B_P\cross_{\tau,\omega}P)$.
Then $V := \overline\sigma\circ i_P$ is a partial isometric
$\omega$\ndash representation of $P$ which satisfies \eqref{eq:V covariant},
and applying $\overline\sigma^{-1}$ we see that $i_P$ is as well.
Since $i_P(s)i_P(s)^* = i_{B_P}(1_s)$ for every $s\in P$,
$i_P$ also satisfies \eqref{eq:V covariant2}, and is hence bicovariant.
Since $\{1_s:s\in P\}$  generates $B_P$ linearly and
$\{i_{B_P}(a)i_P(t): a\in B_P, t\in P\}$ generates
$\Tt_{\cv}(B_P\cross_{\tau,\omega}P)$ as a $C^*$\ndash al\-ge\-bra,
elements of the form
$i_{B_P}(1_s)i_P(t) = i_P(s)i_P(s)^*i_P(t)$
are also generating.
If $V$ is any bicovariant partial isometric $\omega$-representation of $P$,
then by \cite[Proposition~1.3]{lacarae}
there is a representation $\pi_V$ of $B_P$
such that $\pi_V(1_s) = V_sV_s^*$ for every $s\in P$. 
For any $s,t\in P$ the product $V_tV_s = \omega(t,s)V_{ts}$
is a partial isometry; hence by \cite[Lemma~2]{hw} the projections
$V_sV_s^*$ and $V_t^*V_t$ commute, and we deduce that
$\pi_V(a)V_t^*V_t = V_t^*V_t\pi_V(a)$ for every $a\in B_P$ and $t\in P$.
Further,
\begin{align*}
\pi_V(\tau_s(1_t))
& = \pi_V(1_{st})
  = V_{st}V_{st}^* \\
& = (\overline{\omega(s,t)}V_sV_t)(\overline{\omega(s,t)}V_sV_t)^*
  = V_sV_tV_t^*V_s^*
  = V_s\pi_V(1_t)V_s^*,
\end{align*}
so $\pi_V(\tau_s(a)) = V_s\pi_V(a)V_s^*$ for every $s\in P$ and $a\in B_P$.
Thus $(\pi_V,V)$ is a Nica-Toeplitz covariant representation of
$(B_P,P,\tau,\omega)$.
The representation $V_* := \pi_V\times V$ satisfies $V_*\circ i_P = V$.
\end{proof}

We  say that a bicovariant partial isometric $\omega$\ndash representation $V$
is {\em universal\/} if,
for every bicovariant partial isometric $\omega$\ndash representation $W$,
there is a homomorphism of $C^*\{V_s: s\in P\}$
which maps $V_s$ to $W_s$ for each $s\in P$.

\begin{theorem}\label{theorem:bicovariant}
Suppose $(G,P) = *(G_i, P_i)$ is a free product of abelian lattice-ordered groups
and $V$ is a bicovariant partial isometric $\omega$\ndash representation of $P$.
Then $V$ is universal if and only if
\[
\prod_{l=1}^m (V_rV_r^* - V_{rt_l}V_{rt_l}^*)\prod_{k=1}^n (1 - V_{s_k}^*V_{s_k}) \ne 0
\]
whenever $r\in P$, $m,n \ge 1$, and $s_1$, \dots, $s_n$,
$t_1$, \dots, $t_m\in P\setminus\{e\}$.
\end{theorem}

\begin{proof} $V$ is universal if and only if the representation
$V_* = \pi_V\times V$ of $\Tt_{\cv}(B_P\cross_{\tau,\omega}P)$
is faithful.  By Theorem~\ref{theorem:crossed products}, this occurs if and only
if $\pi_V$ acts faithfully on the range of
$\prod_{k=1}^n (1 - V_{s_k}^*V_{s_k})$
whenever $s_1,\dots, s_n\in P\setminus\{e\}$,
and the result follows from \cite[Proposition~1.3]{lacarae}.
\end{proof}

Let $\FF_\infty$ be the free group on infinitely many generators
$z_1$, $z_2$, \dots, and let $\FF_\infty^+$ be the subsemigroup (with identity)
generated by the $z_i$;
the pair $(\FF_\infty, \FF_\infty^+)$ is quasi-lattice ordered.
In \cite{lacarae}, Laca and Raeburn realized the Cuntz
algebra $\Oo_\infty$ as the universal $C^*$\ndash al\-ge\-bra for covariant isometric
representations of $\FF_\infty^+$, and used their characterization
of the faithful representations of $B_P\cross_\tau P$
to derive Cuntz's simplicity result.
We finish by showing that
the universal $C^*$\ndash al\-ge\-bra for bicovariant partial isometric
representations of $\FF_\infty^+$ is reminiscent of $\Oo_\infty$,
and we derive a Cuntz-Krieger-type uniqueness theorem.

First some notation.
For a multi-index $\mu = (\mu_1,\dots, \mu_n)$ we write
$z_\mu := z_{\mu_1}\dotsm z_{\mu_n}$,
and we identify $\FF_\infty^+$ with the set of multi-indices
under concatenation via $z_\mu \leftrightarrow \mu$.

\begin{prop}
Suppose $S$ is a partial isometric representation of $\FF_\infty^+$
in a $C^*$\ndash al\-ge\-bra $B$; that is, $S$ is a semigroup homomorphism
and each $S_\mu$ is a partial isometry.
Then $C^*\{S_\mu: \mu\in\FF_\infty^+\}$ is generated by
$\{S_n: n \in \NN\}$,
and $S$ is bicovariant if and only if

\textup{(a)} the range projections $s_ks_k^*$ for $k\in\NN$
are pairwise orthogonal, and

\textup{(b)} the initial projections $s_k^*s_k$ for $k\in\NN$
are pairwise orthogonal.
\end{prop}

\begin{proof} The first statement is obvious.
In the left-invariant partial order on $\FF_\infty$,
two elements $\mu,\nu\in\FF_\infty^+$ have a common upper bound
if and only if one is an initial word of the other, and then
the least upper  bound is the longer of the two words.
We will show that (a) holds if and only if
\begin{equation}\label{eq:S covariant}
S_\mu S_\mu^* S_\nu S_\nu^*
= \begin{cases}
S_\mu S_\mu^*
  & \text{if $\nu^{-1}\mu \in \FF_\infty^+$,} \\
S_\nu S_\nu^*
  & \text{if $\mu^{-1}\nu \in \FF_\infty^+$,} \\ 
0 & \text{otherwise;} 
\end{cases}
\end{equation}
of course a similar statement holds for (b)
using the right-invariant partial order, and together these prove
the Proposition.

To begin with, \eqref{eq:S covariant} implies (a) since distict generators
of $\FF_\infty^+$ are not comparable.
For the converse, first suppose $\nu^{-1}\mu \in \FF_\infty^+$;
since $S_\nu$ is a partial isometry, we then have
\[
S_\mu S_\mu^* S_\nu S_\nu^*
= S_\mu S_{\nu^{-1}\mu}^* S_\nu^* S_\nu S_\nu^*
= S_\mu S_{\nu^{-1}\mu}^* S_\nu^*
= S_\mu S_\mu^*.
\]
The case $\mu^{-1}\nu \in \FF_\infty^+$ is similar.
Finally, suppose $\mu$ and $\nu$ are not comparable.
Then there exists $\sigma, \mu',\nu'\in\FF_\infty^+$
such that $\mu = \sigma\mu'$,
$\nu = \sigma\nu'$, and $\mu'_1 \ne \nu'_1$.
Condition (a) implies that $S_{\mu'}^* S_{\nu'} = 0$,
and by \cite[Lemma~2]{hw} the range projection of $S_{\nu'}$ commutes with the
initial projection of $S_\sigma$, so
\[
S_\mu^* S_\nu
= S_{\mu'}^* S_\sigma^* S_\sigma S_{\nu'}
= S_{\mu'}^* S_\sigma^* S_\sigma S_{\nu'} S_{\nu'}^* S_{\nu'}
= S_{\mu'}^* S_{\nu'} S_{\nu'}^* S_\sigma^* S_\sigma S_{\nu'}
= 0.
\]

\end{proof}

\begin{theorem}\label{theorem:free}
A bicovariant partial isometric representation $S$
of $\FF_\infty^+$ is universal if and only if each $S_\mu$ is nonzero.
\end{theorem}

\begin{proof}
Suppose each $S_\mu$ is nonzero.
To see that $S$ is universal, we apply  Theorem~\ref{theorem:bicovariant}.
If $\nu\in\FF_\infty^+$, $m,n\ge 1$, and
$\sigma_1$, \dots, $\sigma_m$, $\tau_1$, \dots, $\tau_n\in \FF_\infty^+\setminus\{e\}$,
then we can choose $i,j\in\NN$ such that none of the multi-indices $\sigma_l$
begins with $i$, and none of the multi-indices $\tau_k$ ends with $j$.
Then
\[
\prod_{l=1}^m (S_\nu S_\nu^* - S_{\nu\sigma_l}S_{\nu\sigma_l}^*)
\prod_{k=1}^n (1 - S_{\tau_k}^*S_{\tau_k})
 \ge S_\nu S_iS_i^*S_\nu^* S_j^*S_j
 = S_j^*S_{j\nu i}S_{j\nu i}^*S_j
\]
is nonzero since
$S_j(S_j^*S_{j\nu i}S_{j\nu i}^*S_j)S_j^* = S_{j\nu i}S_{j\nu i}^* \ne 0$.
Hence $S$ is universal.

Now define
$T:\FF_\infty^+\to B(\ell^2(\FF_\infty^+)\otimes\ell^2(\FF_\infty^+))$
by
\[
T_\mu(\delta_\sigma \otimes \delta_\nu)
= \begin{cases} \delta_\sigma \otimes \delta_{\mu\nu}
  & \text{if $\sigma$ ends in $\mu\nu$} \\
  0 & \text{otherwise.} \end{cases}
\]
Then $T$ is a bicovariant partial isometric representation
of $\FF_\infty^+$ in which each $T_\mu$ is nonzero.
If $S$ is universal, then $S_\mu \mapsto T_\mu$ extends
to a homomorphism of $C^*\{S_\mu\}$, and hence each $S_\mu$
must be nonzero.
\end{proof}

\end{document}